\documentclass{amsart}
\usepackage{amssymb}
\newtheorem{tm}{Theorem}
\newtheorem{df}{Definition}
\newtheorem{lm}{Lemma}
\newtheorem{kor}{Corollary}

\newtheorem{conj}{Conjecture}

\theoremstyle{remark}
\newtheorem{remark}{Remark}

\begin{document}
\title[A polynomial variant of a problem of Diophantus]{A POLYNOMIAL VARIANT OF A PROBLEM OF DIOPHANTUS AND ITS CONSEQUENCES}
\author{ALAN FILIPIN AND ANA JURASI\'{C}}

\begin{abstract} We prove that every Diophantine quadruple in $\mathbb{R}[X]$ is regular. More precisely, we prove that if  $\{a, b, c, d\}$ is a set of four non-zero polynomials from $\mathbb{R}[X]$, not all constant, such that the product of any two of its distinct elements increased by $1$ is a square of a polynomial from $\mathbb{R}[X]$, then $$(a+b-c-d)^2=4(ab+1)(cd+1).$$

One consequence of this result is that there does not exist a set
of four non-zero polynomials from $\mathbb{Z}[X]$, not all
constant, such that a product of any two of them increased by a
positive integer $n$, which is not a perfect square, is a square
of a polynomial from $\mathbb{Z}[X]$. Our result also implies that
there does not exist a set of five non-zero polynomials from
$\mathbb{Z}[X]$, not all constant, such that a product of any two
of them increased by a positive integer $n$, which is a perfect
square, is a square of a polynomial from
$\mathbb{Z}[X]$.\end{abstract}

\maketitle


\section{INTRODUCTION}

Diophantus of Alexandria \cite{dio} noted that the product of any
two elements of the set $\big\{\frac{1}{16}, \frac{33}{16},
\frac{17}{4},\frac{105}{16}\big\}$ increased by $1$ is a square of
rational number. A set consisting of $m$ positive integers
(rational numbers) with the property that the product of any two
of its elements increased by 1 is a square of integer (rational
number) is therefore called a Diophantine $m$-tuple. The first
Diophantine quadruple of integers, the set $\{1, 3, 8, 120\}$, was
found by Fermat.

One of the questions of interest is how large those sets can be. No
upper bound for the size of such sets of rational numbers is
known. Gibbs \cite{gib} found some rational Diophantine sextuples.
Very recently, Dujella et al. \cite{duje_rat} proved that there
exist infinitely many rational Diophantine sextuples. In integer
case, which is the most studied, very recently He, Togb\'{e} and
Ziegler \cite{HTZ} announced the proof of the folklore conjecture that there does
not exist a Diophantine quintuple. There is also a stronger version
of that conjecture which states that every Diophantine triple can
be extended to a quadruple with a larger element in a unique way
(see \cite{duje1}):

\begin{conj}\label{conj:1}
If $\{a,b, c, d\}$ is a Diophantine quadruple of integers and $d>\rm{max}\{\textit{a,b,c}\}$,
then $d = d_{+}=a+b+c+2(abc+\sqrt{(ab+1)(ac+1)(bc+1)}).$
\end{conj}

This conjecture is still open. In 1979, Arkin, Hoggatt and Strauss
\cite{AHS} proved that every Diophantine triple of integers
$\{a,b,c\}$ can be extended to a Diophantine quadruple of integers
$\{a,b,c, d_{+}\}$. Baker and Davenport \cite{B-D} proved
Conjecture \ref{conj:1} for the triple $\{a,b,c\}=\{1,3,8\}$ with
the unique extension $d=120$. Many other results are also known
(see \cite{D-P, D-pdeb1,Fujk,B-D-M}) which supports this
conjecture.

Many generalizations of the original problem of Diophantus were
also considered, for example by adding a fixed integer $n$ instead
of $1$, looking at $k$th powers instead of squares, or considering
the problem over domains other than $\mathbb{Z}$ or $\mathbb{Q}$.
We have the following definition:

\begin{df}\label{df1}
Let $m\geq 2$, $k\geq 2$ and let $R$ be a commutative ring with $1$. Let $n\in R$ be a non-zero element and let $\{a_{1}, \ldots ,a_{m}\}$ be a set
of $m$ distinct non-zero elements from $R$ such that
$a_{i}a_{j}+n$ is a $k$th power of an element of $R$ for $1\leq
i<j\leq m$. The set $\{a_{1}, \ldots ,a_{m}\}$ is called a $k$th power Diophantine $m$-tuple with the property $D(n)$ or simply a $k$th power $D(n)$-$m$-tuple in $R$.
\end{df}

It is interesting to find upper bounds for the number of elements of
such sets. Dujella \cite{duje2,duje3} found such bounds for the
integer case and for $k=2$. For other similar results see
\cite{dif1,dif2,dij,gib}. Brown \cite{brown} proved that if
$n$ is an integer, $n\equiv 2\ ({\rm mod}\ {4})$, then
there does not exist a Diophantine quadruple of integers with the
property $D(n)$. Furthermore, Dujella \cite{duje-gen} proved that
if an integer $n$,  $n\not\equiv 2\ ({\rm mod}\ {4})$ and $n\notin S=\{-4,-3,-1,3,5,8,12,20\}$, then
there exists at least one Diophantine quadruple of integers with
the property $D(n)$, and if $n\notin S\cup T$, where
$T=\{-15,-12,-7,7, 13, 15, 21, 24, 28, 32, 48, 60, 84\}$, then
there exist at least two distinct Diophantine quadruples of
integers with the property $D(n)$. For some integers $n$ the
question of the existence of such a Diopantine quadruple is still
unanswered, as it is stated in Dujella's conjecture
\cite{duje-exc}:
\begin{conj}\label{conj:2}
For $n\in S=\{-4,-3,-1,3,5,8,12,20\}$ there does not exist a $D(n)$-quadruple of natural numbers.
\end{conj}

This problem (see \cite[Remark 3]{duje-gen}), asking if there
exist a $D(n)$-quadruple, can be reduced to elements of the set $S'=\{-3,
-1, 3, 5, 8, 20\}$.
\bigskip

In this paper we will consider a polynomial variant of the
problem. A polynomial variant of the problem of Diophantus was
first studied by Jones \cite{jon1,jon2} for the case
$R=\mathbb{Z}[X]$, $k=2$ and $n=1$. There were also considered a
lot of other variants of such a polynomial problem (see
\cite{dif2,dij, dift,difw,dil1,jaialan}). In case of $R$ a
polynomial ring it is usually assumed that, for constant $n$, not
all polynomials in such a $D(n)$-$m$-tuple are constant. In this
paper we first consider the case where $R=\mathbb{R}[X]$, $k=2$
and $n=1$.  Then, we apply the obtained result to the case where
$R=\mathbb{Z}[X]$, $k=2$ and $n$ is a positive integer, to get
other interesting results. We need next two definitions (see \cite{gib-dodano}). Let
$\{a,b,c\}$ be a $D(n)$-triple in a polynomial ring $R$ such that
\begin{eqnarray}\label{prva}
ab+n=r^{2},\
ac+n=s^{2},\ bc+n=t^{2},
\end{eqnarray}
where $r,s,t\in R$ and $n\in\mathbb{Z}\backslash\{0\}$.
\begin{df}\label{def1}
A $D(n)$-triple $\{a,b,c\}$ in $R$ is called regular if
\begin{equation} \label{zv1}
(c-b-a)^{2}=4(ab+n).
\end{equation}
\end{df}Equation \eqref{zv1} is symmetric under permutations of $a$, $b$,
$c$. From \eqref{zv1}, using  \eqref{prva}, we get
\begin{equation}\label{c-reg}c_{\pm}=a+b \pm 2r, \end{equation}
\begin{equation}\label{st-reg}ac_{\pm} + n = (a\pm r)^2, \, bc_{\pm} + n = (b\pm r)^2. \end{equation} \begin{df}\label{def2}
A $D(n)$-quadruple $\{a,b,c,d\}$ in $R$ is called regular if
\begin{eqnarray}\label{drugo}
n(d+c-a-b)^{2}=4(ab+n)(cd+n).
\end{eqnarray}
\end{df}Equation \eqref{drugo} is symmetric under permutations of $a$,
$b$, $c$, $d$. The right hand side of \eqref{drugo} is a square,
so in $\mathbb{Z}[X]$ a regular $D(n)$-quadruple exists only for
$n$ which is a perfect square, whereas regular $D(n)$-triples
exist for every $n$. In $\mathbb{R}[X]$ a regular $D(n)$-quadruple
exists for every positive integer $n$. For example, the set $$\{X,
4X+4\sqrt{5},
9X+6\sqrt{5},\frac{144}{5}X^3+48\sqrt{5}X^2+124X+20\sqrt{5}\}$$ is
a regular $D(5)$-quadruple in $\mathbb{Q}(\sqrt{5})[X]$.  In
$\mathbb{C}[X]$ a regular $D(n)$-quadruple exists for every
nonzero integer $n$. Equation (\ref{drugo}) is a quadratic
equation in $d$ with roots
\begin{equation}\label{d-plus}d_{\pm}=a+b+c+\frac{2}{n}(abc\pm rst),\end{equation}
and it holds
\begin{equation}\label{d-plus-minus}ad_{\pm}+n=\frac{1}{n}u_{\pm}^2, \, bd_{\pm}+n=\frac{1}{n}v_{\pm}^2, \,
cd_{\pm}+n=\frac{1}{n}w_{\pm}^2,\end{equation}where \begin{eqnarray}\label{uvw}u_{\pm}=at \pm rs,\  v_{\pm}=bs \pm rt,\ w_{\pm}=cr \pm st.\end{eqnarray}
An irregular $D(n)$-quadruple in $R$ is one that is not
regular.  It is known from \cite{AHS} (all relations are obtained using only algebraic manipulations so they hold in every ring $R$) that every $D(1)$-pair $\{a,b\}$ in $R$ can be extended to a regular $D(1)$-quadruple in $R$: \begin{eqnarray}\label{pr}\{a,b,c_{+},d_{+}\},\end{eqnarray}where $d_{+}=4r(a+r)(b+r)$, $d_{-}=0$. Notice that if $d_{-}=0$, then from (\ref{drugo}) it follows that the $D(1)$-triple $\{a,b,c\}$ is regular.

For the simplicity, in the rest of the paper we use a term polynomial $D(1)$-$m$-tuple
for a second power $D(1)$-$m$-tuple in $\mathbb{R}[X]$. We also
use a term $D(n)$-$m$-tuple in $\mathbb{Z}[X]$ for the second power $D(n)$-$m$-tuple in $\mathbb{Z}[X]$, and we will always specify that it is from $\mathbb{Z}[X]$. From \cite{dij} it
follows that there are at most $7$ elements in a polynomial
$D(1)$-$m$-tuple and also in a $D(n)$-$m$-tuple from
$\mathbb{Z}[X]$ for $n\in \mathbb{Z}\setminus\{0\}$. Dujella and
Fuchs \cite{dif2} proved that every $D(1)$-quadruple in
$\mathbb{Z}[X]$ is regular, i.e. there are at most $4$ elements in
a $D(1)$-$m$-tuple from $\mathbb{Z}[X]$.  We furthermore extend their result:
\begin{tm}\label{Tm1}
Every polynomial $D(1)$-quadruple is regular.\end{tm}

Suppose that $\{a,b,c,d\}$ is an irregular $D(n)$-quadruple in
$\mathbb{Z}[X]$ for $n\in \mathbb{N}$. Then, the set
$\{\frac{a}{\sqrt{n}},\frac{b}{\sqrt{n}},\frac{c}{\sqrt{n}},\frac{d}{\sqrt{n}}\}$
would be an irregular polynomial $D(1)$-quadruple, which is not
possible by Theorem \ref{Tm1}. For every $n\in \mathbb{N}$ which
is a perfect square a regular $D(n)$-quadruple in $\mathbb{Z}[X]$
can be obtained. For example, from $D(1)$-quadruple (\ref{pr}) by
multiplying its elements with $\sqrt{n}$. Since for $n\in
\mathbb{N}$ which is not a perfect square there does not exist a
regular $D(n)$-quadruple in $\mathbb{Z}[X]$, from Theorem
\ref{Tm1} we have:

\begin{kor}\label{Kor1}
There does not exist a $D(n)$-quadruple in $\mathbb{Z}[X]$ for a
positive integer $n$ which is not a perfect square. Furthermore,
there does not exist a $D(n)$-quintuple in $\mathbb{Z}[X]$ for a
positive integer $n$ which is a perfect square.\end{kor}

Let us mention that there exist a polynomial $D(n)$-sextuple in $\mathbb{Z}[X]$, where $n$ is not a constant polynomial (see \cite{duj_kaz,duje_rat}). Moreover, in all those examples $n$ is a square
in $\mathbb{Z}[X]$, while for $n\in\mathbb{Z}[X]$ non-square, there exist examples of $D(n)$-quintuples in $\mathbb{Z}[X]$ (see \cite{dif3}).

Dujella and Fuchs \cite{dif1} proved that there are at most $3$
elements in a $D(-1)$-$m$-tuple in $\mathbb{Z}[X]$, i.e. they
proved a polynomial variant of Conjecture \ref{conj:2} for $n=-1$.
By Corollary \ref{Kor1} we prove a polynomial variant of
Conjecture \ref{conj:2} for $n\in \{3,5,8,12,20\}$. For an integer $n<0$ we can not apply Theorem \ref{Tm1} to observe a polynomial $D(n)$-quadruple in $\mathbb{Z}[X]$ because in that case a $D(1)$-quadruple $\{\frac{a}{\sqrt{n}},\frac{b}{\sqrt{n}},\frac{c}{\sqrt{n}},\frac{d}{\sqrt{n}}\}$ is from $\mathbb{C}[X]$. Notice that in
a polynomial variant of Conjecture \ref{conj:2} we can not reduce
the set $S$ to the set $S'$.

In order to prove Theorem \ref{Tm1} we consider the ring
$\mathbb{R}[X]$, using the relation "$<$" between its elements. We
partially follow the strategy used in \cite{dif2} for
$\mathbb{Z}[X]$. However, not everything is the same in $\mathbb{R}[X]$, so we need to introduce some new ideas. In Section 2 we transform the problem of
extending a polynomial $D(1)$-triple $\{a,b,c\}$ to a polynomial
$D(1)$-quadruple $\{a,b,c,d\}$ into solving a system of
simultaneous Pellian equations, which reduces to finding
intersections of binary recurring sequences of polynomials. The
main difference from approach in \cite{dif2} is that when
minimality is observed we have to consider relation "$\leq$"
between the degrees of polynomials instead of between polynomials.
On the other hand, we use some results from \cite{dij} valid for a
polynomial $D(1)$-$m$-tuple in $\mathbb{C}[X]$, where we do not
have the relation "$<$" between its elements. In Section 3 we find
a gap principle for degrees of elements in a polynomial
$D(1)$-triple $\{a,b,c\}$ and we also describe all possible
initial terms of the recurring sequences obtained for such triple.
Using results from Sections 2 and 3, in Section 4 we prove Theorem
\ref{Tm1}.

\section{REDUCTION TO INTERSECTIONS OF RECURSIVE SEQUENCES}

Let $\mathbb{R}^{+}[X]$ denote the set of all polynomials with
real coefficients with positive leading coefficient. For
$a,b\in\mathbb{R}[X]$, $a<b$ means that $b-a\in\mathbb{R}^{+}[X]$. For $a\in\mathbb{R}[X]$, we define $|a|=a$ if $a\geq 0$ and $|a|=-a$ if $a< 0$.

Let us consider an arbitray extension of a polynomial $D(1)$-triple $\{a,b,c\}$, where $0<a<b<c$, to a polynomial $D(1)$-quadruple $\{a,b,c,d\}$. We first observe equations (\ref{prva}) for $n=1$ and $r,s,t\in\mathbb{R}^{+}[X]$. Let $A, B, C, R, S, T$ be the leading coefficients of polynomials $a, b, c, r, s, t$, respectively. By (\ref{prva}), $AB=R^{2}$, $AC=S^{2}$ and $BC=T^{2}$.  Hence, there is no loss of generality in assuming that $a,b,c\in\mathbb{R}^{+}[X]$.  By \cite[Lemma 1]{dil1}, there is at most one constant in a polynomial $D(1)$-$m$-tuple in $\mathbb{C}[X]$. In the proof it is used the famous theorem of Mason \cite{mas, sto} (see also \cite{dil1}), usually called the \textit{abc} theorem for polynomials. Since any polynomial $D(1)$-$m$-tuple in $\mathbb{R}[X]$ is also a $D(1)$-$m$-tuple in $\mathbb{C}[X]$, it follows that a polynomial $D(1)$-$m$-tuple can not contain two constants. Let us denote by $\alpha,\beta,\gamma$ degrees of $a,b,c$, respectively. Hence, $0\leq\alpha\leq\beta\leq\gamma$ and $\beta,\gamma>0$. Let \begin{eqnarray}\label{treca}ad+1=x^{2},\ \  bd+1=y^{2},\ \ cd+1=z^{2},\end{eqnarray}where $x,y,z\in\mathbb{R}[X]$. Note that $x$, $y$ and $z$ can be $<0$, because by taking only positive values we would exclude some possibilities obtained for polynomials from $\mathbb{C}[X]$. By (\ref{prva}) and (\ref{treca}), $\alpha,\beta,\gamma,\delta$ are all even or all odd numbers and $d\in\mathbb{R}^{+}[X]$.

Eliminating $d$ from (\ref{treca}), we obtain the system of Pellian equations\begin{eqnarray}\label{cetvrta}az^{2}-cx^{2}&=&a-c,\\
bz^{2}-cy^{2}&=&b-c.\label{peta}
\end{eqnarray}We want to find solutions $(z,x)$ and $(z,y)$ of (\ref{cetvrta}) and (\ref{peta}), respectively. The following lemma describes these solutions.

\begin{lm}\label{Lm1}Let $(z,x)$ and $(z,y)$ be solutions, with $x,y,z\in\mathbb{R}[X]$, of  (\ref{cetvrta}) and (\ref{peta}), respectively. Then there exist solutions $(z_0,x_0)$ and $(z_1,y_1)$, with $z_0, x_0, z_1, y_1\in\mathbb{R}[X]$, of (\ref{cetvrta}) and (\ref{peta}), respectively, such that:
\begin{eqnarray}\label{ograda1}|x_0|\geq 1, &\ \ & |z_0|\geq 1,\\
|y_1|\geq 1, &\ \ \label{ograda2}& |z_1|\geq 1\end{eqnarray}and \begin{eqnarray}\label{ograda3}{\rm deg}(z_{0})\leq \frac{3\gamma -\alpha}{4}, &\ \ & {\rm deg}(x_{0})\leq\frac{\alpha +\gamma}{4},\\
{\rm deg}(z_{1})\leq \frac{3\gamma -\beta}{4}, &\ \ \label{ograda4}& {\rm deg}(y_{1})\leq\frac{\beta +\gamma}{4}.\end{eqnarray}There also exist non-negative integers $m$ and $n$ such that
\begin{equation}\label{deseto}
 z\sqrt{a}+x\sqrt{c}=(z_{0}\sqrt{a}+x_{0}\sqrt{c})(s+\sqrt{ac})^{m},
 \end{equation}
\begin{equation}\label{jedanaesto}
 z\sqrt{b}+y\sqrt{c}=(z_{1}\sqrt{b}+y_{1}\sqrt{c})(t+\sqrt{bc})^{n}.
 \end{equation}
\end{lm}
\emph{Proof.}
The statements (\ref{ograda3}) - (\ref{jedanaesto}) follow directly from \cite[Lemma 4]{dil1}.

From the proof of \cite[Lemma 4(v)]{dil1}, we have that if $c|(z^2-1)$, then $c|(z_0^2-1)$. Hence, there exists
$d_{0}\in\mathbb{R}[X]$ such that $cd_0=z_0^2-1$. Then, by (\ref{cetvrta}), $ad_0=x_0^2-1$. Therefore, \begin{equation}\label{d_0}ad_{0}+1=x_{0}^{2} \ \ \rm{and}\ \  \textit{cd}_{0}+1=\textit{z}_{0}^{2}.\end{equation}
Since $c>0$ and for $d_0\neq 0$ we have $\rm {deg}(\textit{z}_{0})>0$,  by (\ref{d_0}), we conclude that $d_{0}\geq 0$. Further,  by (\ref{d_0}), $x_0^2\geq 1$ and $z_0^2\geq 1$ so (\ref{ograda1}) holds. The proof of the statements of the lemma for $z_1$ and $y_1$ are analogous. \hfill$\square$\\

In particular, if $z_0$ is a constant then $z_0=\pm 1$. By \cite[Lemma 4]{dil1}, if $(z_0,x_0) \neq  (\pm 1, \pm 1)$, then $\rm {deg}(\textit{z}_{0})\geq\frac{\gamma}{2}$ and $\rm {deg}(\textit{x}_{0})\geq\frac{\alpha}{2}$. Analogously, if $(z_1,y_1) \neq  (\pm 1, \pm 1)$, then $\rm {deg}(\textit{z}_{1})\geq\frac{\gamma}{2}$ and $\rm {deg}(\textit{y}_{1})\geq\frac{\beta}{2}$.

By Lemma \ref{Lm1},  $z=v_{m}=w_{n}$, where the sequences
$(v_{m})$ and $(w_{n})$ are, for $m,n\geq 0$, defined by
\begin{eqnarray}\label{osma}
v_{0}&=&z_{0},\ \ v_{1}=sz_{0}+cx_{0},\ \ v_{m+2}=2sv_{m+1}-v_{m},\\
w_{0}&=&z_{1},\ \ w_{1}=tz_{1}+cy_{1},\ \ \label{deveta}
w_{n+2}=2tw_{n+1}-w_{n}.\end{eqnarray}Initial values $(z_0,x_0)$
and $(z_1,y_1)$ from (\ref{osma}) and (\ref{deveta}) are some
solutions of (\ref{cetvrta}) and (\ref{peta}), respectively, with estimates
(\ref{ograda3}) and (\ref{ograda4}) satisfied.

\begin{remark}\label{rem2}If we observe the equation $v_{m}=w_{n}$ just for $z_{0}>0$ or just for $z_{0}<0$ we lose some solutions of Pellian equation (\ref{cetvrta}). The analogous situation is for $z_{1}$.  Without loss of generality we may assume that $x_0>0$ because, by (\ref{treca}), for $z$ and $-z$ we obtain the same $d$. The analogue situation is for $y_1$.
\end{remark}

By \cite[Lemma 5]{dil1}, it
follows
that
\begin{equation}\label{deseta}
{\rm deg}(v_{m})=(m-1)\frac{\alpha
+\gamma}{2}+{\rm deg}(v_{1}),
\end{equation}
for $m\geq 1$
and
\begin{equation}\label{jedanaesta}
\frac{\gamma}{2}\leq {\rm deg}(v_{1})\leq\frac{\alpha
+5\gamma}{4}.
\end{equation}
Similarly,
\begin{equation}\label{dvanaesta}
{\rm deg}(w_{n})=(n-1)\frac{\beta
+\gamma}{2}+{\rm deg}(w_{1}),
\end{equation}
for $n\geq 1$
and
\begin{equation}\label{trinaesta}
\frac{\gamma}{2}\leq {\rm deg}(w_{1})\leq\frac{\beta
+5\gamma}{4}.
\end{equation}

For the rest of the paper we assume that $\{a,b,c,d'\}$ is an irregular polynomial $D(1)$-quadruple with $d'>c$ and we try to prove that such quadruple does not exist. Let us denote by $\delta$ degree of $d'$. By \cite[Lemma 5]{dij}, we know that \begin{equation}\label{ograda-delta}\delta \geq
\displaystyle\frac{3\beta +5\gamma}{2},\end{equation} since
$\mathcal{D}_{p}=\big\{\frac{\sqrt{-3}}{2}, -\frac{2\sqrt{-3}}{3}(p^{2}-1),
\frac{-3+\sqrt{-3}}{3}p^{2}+\frac{2\sqrt{-3}}{3},\frac{3+\sqrt{-3}}{3}p^{2}+\frac{2\sqrt{-3}}{3}
\big\}$, where $p\in\mathbb{C}[X]$ is a
non-constant polynomial, is not from $\mathbb{R}[X]$. Assume that $\delta$ is minimal possible degree for which (\ref{ograda-delta}) holds.

\section{GAP PRINCIPLE FOR DEGREES AND PRECISE DETERMINATION OF INITIAL TERMS}

In this
section we describe in a gap principle all possible relations
between $\alpha$, $\beta$ and $\gamma$.  In Lemma \ref{Lm1} we
have some useful facts about initial terms of the recurring
sequences $(v_{m})$ and $(w_{n})$. Furthermore, in this section we determine
all possible initial terms of the sequences, for the triple
$\{a,b,c\}$ described in Section 2. We need  the following expressions proved for polynomials from $\mathbb{Z}[X]$, but the constructions also works in $\mathbb{R}[X]$. From \cite[Lemma 3]{dift} and (\ref{uvw}), we have $u_{-},v_{-}<0$. Also, $u_{+},v_{+}>0$.  From \cite[Lemma 1]{dift}, \begin{equation}c=a+b+ d_{\pm}+2(ab d_{\pm}\mp ru_{\pm}v_{\pm}).\label{ccc}\end{equation}By (\ref{ccc}), we have $c=e_\mp$, where $e_\mp$ is obtained by applying (\ref{d-plus}) on the $D(1)$-triple $\{a,b,d_\pm\}$. For $d_{-}=0$, by (\ref{ccc}), we obtain $c=c_+$. Otherwise, by (\ref{d-plus-minus}), for $n=1$, it holds $d_{-}> 0$.  By (\ref{ccc}) with the lower signs, we conclude that $c>a+b$ so $c^2>c(a+b)+1$. By that and (\ref{prva}), $s^2t^2=abc^2+c(a+b)+1<c^2r^2$. Therefore, $w_{-}=cr-st>0$. Also, $w_{+}>0$. From (\ref{ccc}), using (\ref{prva}) and (\ref{uvw}), we get\begin{equation}c=a+b- d_{\pm}+2rw_{\pm}.\label{w}\end{equation}

For $d_{-}\neq 0$, from \cite[Lemma 1]{dfl} and \cite[Lemma 2]{dij}, we have $c>2abd_{-}$, so $0\leq\rm{deg}(\textit{d}_{-})\leq\gamma-\alpha-\beta<\gamma$, i.e. $\gamma\geq\alpha+\beta$. But, we prove even more:

\begin{lm}\label{Lm2}Let $\{a,b,c\}$ be an arbitrary polynomial $D(1)$-triple, where $a<b<c$, and let $d_{-}$ be defined by (\ref{d-plus}) for $n=1$. Then $d_{-}=0$ or $\rm{deg}(\textit{d}_{-})=\gamma-\alpha-\beta.$\end{lm}
\emph{Proof.}
Let $d_{-}\neq 0$. By (\ref{prva}) and (\ref{d-plus-minus}), $\rm{deg}(\textit{abd}_{-})=\alpha+\beta+\rm{deg}(\textit{d}_{-})=\rm{deg}(\textit{ru}_{-}\textit{v}_{-})$. If $\beta<\gamma$, by (\ref{ccc}), $\rm{deg}(\textit{d}_{-})=\gamma-\alpha-\beta$. Let $\beta=\gamma$. Since $\rm{deg}(\textit{abd}_{-})\geq \gamma$, $abd_->0$ and $ru_-v_->0$ then, by (\ref{ccc}), $\rm{deg}(\textit{abd}_{-})=\gamma$ and $\alpha=\rm{deg}(\textit{d}_{-})=0$.\hfill$\square$

\begin{remark}\label{rem3}By the proof of Lemma \ref{Lm2}, if $\beta=\gamma$ then $d_{-}=0$, or if $d_{-}\neq 0$, $\alpha=\rm{deg}(\textit{d}_{-})=0$. The equation $a^{2}+1=u_{-}^2$ has a solution $u_{-}\in\mathbb{R}$ for every $a\in\mathbb{R}$. Since we cannot have two different constants in a polynomial $D(1)$-quadruple, we conclude that if $a$ and $d_{-}$ are non-zero constants, then $d_{-}=a$. This is not possible in $\mathbb{Z}$ for $a\neq 0$. In $\mathbb{Q}$ it is possible for some values of $a$. \end{remark}

Any $D(1)$-triple $\{a,b,c\}$ can obviously be extended to a $D(1)$-quadruple $\{0,a,b,c\}$ and also to a $D(1)$-quadruple $\{a,a,b,c\}$ if $a$ is a constant. A $D(1)$-quadruple with a relaxed condition that its elements need not be distinct and need not be non-zero is called improper\footnote{Note that in an (improper) $D(1)$-quadruple we cannot have two equal non-constant polynomials because then it would be for example $(b-w_{-})(b+w_{-})=-1$ and it is not possible that both factors on the left hand side of this equation are constant.} $D(1)$-quadruple. Such a quadruple can be regular or irregular. Also, any $D(1)$-triple $\{a,b,c\}$ can be extended to regular $D(1)$-quadruples  $\{a,b,c, d_{\pm}\}$. Hence, the equation $v_{m}=w_{n}$ has nontrivial solutions. Moreover, it holds:

\begin{lm}\label{Lm3}Let $\{a,b,c\}$, with $a<b<c$, be a subtriple of an irregular polynomial $D(1)$-quadruple $\{a,b,c,d'\}$, where $d'>c$ with minimal $\delta$. Let $v_m=w_n$ and define $d=\frac{v_ {m}^2-1}{c}$.
\begin{enumerate}
    \item [\textbf{a)}] If $d=d_{-}$, then $v_{m}=w_{n}=\pm w_{-}$ for  $m,n\in\{0,1\}$.
\item [\textbf{b)}] If $d=d'$, then $v_m=w_n=\pm z$ for $m\geq 3$ and $n\geq 3$.
\end{enumerate}
 If $0\in\{m,n\}$, then $d=d_{-}$ or the quadruple $\{a,b,c,d\}$ is irregular with $d=0$ or $d=a$ and $a$ is a constant. If $(m,n)=(1,1)$ then there are four possibilities: $d=d_{-}$; if $\gamma\geq\alpha+2\beta$, it can be $d=d_{+}$; if the quadruple $\{a,b,c,d\}$ is irregular we can have $d=0$ or $d=a$ and $a$ is a constant.
\end{lm}
\emph{Proof.}
\textbf{a)} By Lemma \ref{Lm2} and (\ref{d-plus-minus}), if $d_{-}\neq 0$, then \begin{equation}\label{stupanj_w}\rm{deg}(\textit{w}_{-})=\gamma-\frac{\alpha+\beta}{2}<\gamma.\end{equation}By (\ref{deseta}) and (\ref{dvanaesta}), $\rm{deg}(\textit{v}_{\textit{m}})\geq \gamma$ for $m\geq 2$ and $\rm{deg}(\textit{w}_{\textit{n}})\geq \gamma$ for $n\geq 2$, so $d_{-}$ must arise from $v_{m}=w_{n}$ for $m,n\in\{0,1\}$. If $d_{-}=0$ then $w_{-}=1$, which implies $\rm{deg}(\textit{w}_{-})=0$.

\textbf{b)} If $d=d'$, by \cite[Proposition 1]{dij}, we have $m\geq 3$ and $n\geq 3$. Note that in \cite{dij} the authors considered polynomial $D(1)$-quadruples in $\mathbb{C}[X]$, but every $D(1)$-quadruple in $\mathbb{R}[X]$ is also a $D(1)$-quadruple in $\mathbb{C}[X]$.

If $0\in\{m,n\}$ then, by the proof of \cite[Proposition 1]{dij}, we have $\rm{deg}(\textit{d})<\gamma$. By Lemma \ref{Lm2}, we can have  $d=d_{-}$. By (\ref{d-plus}), we conclude \begin{equation}\label{stupanj_d_plus}\rm{deg}(\textit{d}_{+})=\alpha+\beta+\gamma>\gamma\end{equation} and by (\ref{ograda-delta}), $\rm{deg}(\textit{d}')>\gamma$. By minimality assumption, only possible irregular quadruples are those with $d=0$ or $d=a$ if $a$ is a constant.

Let $v_{1}=w_{1}$. By \cite[Proposition 1]{dij}, $\rm{deg}(\textit{d})<\gamma$ and $\{a,b,c,d\}$ is irregular quadruple (by minimality assumption, then $d=0$ or $d=a$ if $a$ is a constant), or $d=d_{\pm}$. By  and (\ref{d-plus-minus}) and (\ref{stupanj_d_plus}), $\rm{deg}(\textit{w}_{+})=\gamma+\frac{\alpha+\beta}{2}$ so, by  (\ref{jedanaesta}), we can have $d=d_{+}$ and $w_{+}=v_1$ only if $\gamma\geq\alpha+2\beta$.  \hfill$\square$\\

We look for all possible initial terms of the recurring sequences
$(v_{m})$ and $(w_{n})$ for the triple $\{a,b,c\}$, which is a
subtriple of the observed irregular polynomial $D(1)$-quadruple
$\{a,b,c,d'\}$.  We consider $d_{0},d_{1}\in\mathbb{R}[X]$,
where there hold (\ref{d_0}) and
\begin{equation}\label{d_1}bd_{1}+1=y_{1}^{2}, \ \
\textit{cd}_{1}+1=\textit{z}_{1}^{2}.\end{equation} From
(\ref{d_1}), we see that $d_{1}\geq 0$. As we concluded for
$d_{-}$ in Remark \ref{rem3}, if $a$ and $d_{0}$ are non-zero
constants\footnote{This situation is also described in \cite[Lemma 4]{dij}. In $\mathbb{C}[X]$ there are some possibilities which does not exist in $\mathbb{R}[X]$. For example,  $x_0=0$, $a=\pm i$ and $z_0^2=\pm ci+1$.} then $d_{0}=a$. If $\alpha>0$, then $d_0$ and $d_{1}$ can
be constants. Furthermore, obviously we can have $d_{0}=b$ or $d_{0}=d_{-}$ and
$d_{1}=a$ or $d_{1}=d_{-}$.
\begin{remark}\label{rem4}
If $d_0=d_1=d$ then, by (\ref{ograda3}), (\ref{ograda4}), (\ref{stupanj_d_plus}),  (\ref{ograda-delta}) and the minimality assumption, $d=d_{-}$ or  $d=0\neq d_{-}$ or $d=a\neq d_{-}$, where $a$ is a constant.\end{remark}

In the following lemma we consider all possibilities for $d_{-}$. Similar gap principle is well known in classical case and was also used in considering a polynomial variants of the problem of Diophantus (see e.g. \cite[Lemma 4]{dif2}), but here we obtained more information about possible triples. Notice that conclusions about degrees hold for every triple $\{a,b,c\}$, but conclusions about initial terms hold only for the case where $d_{-}$ arise from the intersections of binary recursive sequences.

\begin{lm}\label{Lm4} Let $\{a,b,c\}$ be a polynomial $D(1)$-triple, with $a<b<c$,  for which (\ref{prva}) holds for $n=1$. Then:
\begin{enumerate}
    \item [\textbf{1.}] If $d_{-}=0$, then $d_{0}=d_{1}=0$. In this case $c=a+b+ 2r$ and $\beta=\gamma$. Also, if $\alpha<\beta$, then $C=B$, and if $\alpha=\beta$, then $C=A+B+ 2\sqrt{AB}$.

 \item [\textbf{2.}]
\begin{enumerate}
    \item [\textbf{a)}] If $d_{-}=a\in\mathbb{R}\setminus\{0\}$, then $d_{0}=d_{1}=a\in\mathbb{R}\setminus\{0\}$. In this case $\alpha=0$, $\beta=\gamma$, $c=b+ 2rs$.
 \item [\textbf{b)}] If $d_{-}\in\mathbb{R}\setminus\{0,a\}$, then $d_{0}=d_{1}=d_{-}\in\mathbb{R}\setminus\{0,a\}$. In this case $\alpha>0$ and $\gamma=\alpha+\beta$.
\end{enumerate}

  \item [\textbf{3.}]If $\rm{deg}(\textit{d}_{-})>0$, then we have the following possibilities:
\begin{enumerate}
 \item [\textbf{a)}] $d_{0}=d_{1}=d_{-}$, where $\rm{deg}(\textit{d}_{-})\leq\alpha$,  $\alpha>0$ and $\alpha+\beta<\gamma\leq 2\alpha+\beta$,
 \item [\textbf{b)}] $d_{0}=d_{-}$ and $d_{1}=a$, where $\alpha\leq\rm{deg}(\textit{d}_{-})\leq\beta$,  $\alpha\geq 0$ and $2\alpha+\beta\leq \gamma\leq \alpha+2\beta,$
\item [\textbf{ c)}] $d_{0}=b$ and $d_{1}=d_{-}$, where $\rm{deg}(\textit{d}_{-})=\alpha$,  $\alpha=\beta$ and $\gamma=3 \alpha$,
\item [\textbf{d)}] $d_{0}=b$ and $d_{1}=a$,  where $\beta\leq\rm{deg}(\textit{d}_{-})<\gamma$, $\alpha\geq 0$ and $\gamma\geq \alpha+2\beta.$
\end{enumerate}

\end{enumerate}
\end{lm}
\emph{Proof.} By Lemma \ref{Lm2}, we have $d_{-}=0$ or $\rm{deg}(\textit{d}_{-})=\gamma-\alpha-\beta\geq 0$.\\

\textbf{1.)} We noticed that if $d_{-}=0$, then $c=c_+$, i.e. the triple $\{a,b,c\}$ is regular\footnote{For example, $D(1)$-triples $\{1, X^2+2X,X^2+4X+3\}$ and $\{X-1,X+1,4X\}$.}. Also, $\gamma\leq\beta$, thus $\gamma=\beta$. By (\ref{c-reg}) and (\ref{prva}), if $\alpha<\beta$ then $C=B$, and if $\alpha=\beta$ then $C=A+B+ 2\sqrt{AB}$.

By (\ref{d-plus-minus}), $w_{-}=1$. By (\ref{jedanaesta}) and (\ref{trinaesta}), $\rm{deg}(\textit{v}_{1}), \rm{deg}(\textit{w}_{1}) \geq \frac{\gamma}{2}$. Therefore, $v_{0}=w_{0}=1$. By (\ref{osma}), (\ref{deveta}) and Remark \ref{rem2}, we get $z_{0}=z_{1}=\pm 1$. By (\ref{d_0}) and (\ref{d_1}), $d_{0}=d_{1}=0$.

\textbf{2.)} If $\rm{deg}(\textit{d}_{-})= 0$, then by Lemma \ref{Lm2}, we have
\begin{eqnarray}\label{osma2}\gamma= \alpha+\beta.\end{eqnarray}
\textbf{2. a)} Let $d_{-}=a$ and $\alpha =0$. By (\ref{osma2}), $\beta=\gamma$. By (\ref{d-plus-minus}), we have\footnote{An example of such a case is a $D(1)$-triple $\{\frac{4}{3},\frac{4X^2+2X-2}{3},12X^2+6X\}$, where $d_{-}=\frac{4}{3}$.} $w_{-}= s$, so by (\ref{w}), \begin{equation}\label{a}c=b+2rs.\end{equation}

By Lemma \ref{Lm3} and Remark \ref{rem2}, we consider the cases $v_{m}=w_{n}=\pm s$ for  $m,n\in\{0,1\}$. For $v_{0}=w_{0}=\pm s$, by  (\ref{osma}) and (\ref{deveta}),  $z_{0}=z_{1}=\pm s$.  By (\ref{d_0}) and (\ref{d_1}), $d_{0}=d_{1}=a$. For $v_{0}=w_{1}=\pm s$, by  (\ref{osma}) and (\ref{deveta}),   \begin{equation}\label{z0}z_{0}=tz_{1}+cy_{1}=\pm s.\end{equation} From (\ref{uvw}) and (\ref{z0}), we have $\pm cr\mp st=tz_{1}+cy_{1}$, so \begin{equation}\label{djeljivost}c(\pm r-y_{1})=t(z_{1}\pm s).\end{equation} By (\ref{prva}), $t|(\pm r-y_{1})$. From considering degrees of polynomials, we conclude that $y_{1}=\pm r$ and $z_{1}=\mp s$. By (\ref{d_0}) and (\ref{d_1}), $d_{0}=d_{1}=a$. For $v_{1}=w_{0}=\pm s$, by (\ref{osma}) and (\ref{deveta}),  $sz_{0}+cx_{0}=z_{1}=\pm s$. Hence, $cx_{0}=s(\pm 1-z_{0})$. By (\ref{prva}), we have $s|x_{0}$. Since $x_0\neq 0$, it is not possible, because of the degrees of these polynomials. For $v_{1}=w_{1}=\pm s$, by  (\ref{osma}) and (\ref{deveta}), we have $sz_{0}+cx_{0}=tz_{1}+cy_{1}=\pm s$. As we concluded, this is not possible.

\textbf{2. b)} Let us consider the case\footnote{Such an example is a $D(1)$-triple $\{125X^2+50X, 12500000000X^{10}+26000000000X^9+23070000000X^8+11392000000X^7+3424950000X^6+644520000X^5+75187000X^4+5200000X^3+194525X^2+3300X+15, 1250000000000X^{12}+3100000000000X^{11}+3372000000000X^{10}+2114000000000X^9+844190000000X^8+224024000000X^7+40005200000X^6+4764480000X^5+367264500X^4+17315400X^3+452640X^2+5500X+\frac{96}{5}\}$, where $d_{-}=\frac{1}{5}$.} where $d_{-}\in\mathbb{R}\setminus\{0,a\}$. It holds (\ref{osma2}), and since we can not have two different constants in a $D(1)$-quadruple, $\alpha >0$. By (\ref{uvw}), $w_{-}=cr- st\neq  s$. By Lemma \ref{Lm3} and Remark \ref{rem2}, $v_{m}=w_{n}=\pm cr\mp st$ for  $m,n\in\{0,1\}$. For $v_{0}=w_{0}=\pm cr\mp st$, by  (\ref{osma}) and (\ref{deveta}), we have $z_{0}=z_{1}=\pm cr\mp st$. By (\ref{d_0}) and (\ref{d_1}),  $d_{0}=d_{1}=d_{-}\in\mathbb{R}\setminus\{0,a\}$. For $v_{0}=w_{1}=\pm cr\mp st$, by  (\ref{osma}) and (\ref{deveta}),  \begin{equation}\label{z00}z_{0}=tz_{1}+cy_{1}=\pm cr\mp st.\end{equation}  We obtain (\ref{djeljivost}), so again $y_{1}=\pm r$ and $z_{1}=\mp s$. Hence, $d_1=a$. Using (\ref{d_1}) and (\ref{osma2}), from (\ref{ograda4}) we obtain a contradiction. Analogously as (\ref{z00}), the case $v_{1}=w_{0}=\pm cr\mp st$ is not possible. For $v_{1}=w_{1}=\pm cr\mp st$ we obtain a contradiction analogously as in the previous cases.

\textbf{3.} If $\rm{deg}(\textit{d}_{-})>0$, then by Lemma \ref{Lm2}, we have
\begin{eqnarray}\label{osma3}\gamma>\alpha+\beta.\end{eqnarray}By Lemma \ref{Lm3} and Remark \ref{rem2}, $v_{m}=w_{n}=\pm cr\mp st$ for  $m,n\in\{0,1\}$.

\textbf{3. a)} For the case\footnote{For example, a $D(1)$-triple $\{16X^3-4X, 64X^5-48X^3+8X, 4096X^9+4096X^8-4096X^7-4096X^6+1408X^5+1280X^4-192X^3-128X^2+9X+3\}$, where $d_{-}=X+1$.} $v_{0}=w_{0}=\pm cr\mp st$, by  (\ref{osma}) and (\ref{deveta}),  $z_{0}=z_{1}=\pm cr\mp st$. By (\ref{d_0}) and (\ref{d_1}), $d_{0}=d_{1}=d_{-}$. Using (\ref{stupanj_w}) and (\ref{ograda4}), we obtain $\gamma\leq 2\alpha+\beta$. Then, by Lemma \ref{Lm2}, $\rm{deg}(\textit{d}_{-})\leq\alpha.$ Therefore, $\alpha >0$.

\textbf{3. b)}  For the case\footnote{For example, the $D(1)$-triple $\{\frac{1}{5},12500000000X^{10}+26000000000X^9+23070000000X^8+11392000000X^7+3424950000X^6+644520000X^5+75187000X^4+5200000X^3+    194525X^2+ 3300X+15,1250000000000X^{12}+3100000000000X^{11}+3372000000000X^{10}+2114000000000X^9+844190000000X^8+224024000000X^7+40005200000X^6+4764480000X^5+367264500X^4+17315400X^3+452640X^2+5500X+\frac{96}{5}\}$, where $d_{-}=125X^2+50X$.} $v_{0}=w_{1}=\pm cr\mp st$, by  (\ref{osma}) and (\ref{deveta}),  we have (\ref{z00}). Hence, $d_{0}=d_{-}$. (For $\alpha=0$ and $d_{-}=a$ we have the case 2.a).)  By (\ref{stupanj_w}) and (\ref{ograda3}), $\gamma\leq \alpha+2\beta$. By Lemma \ref{Lm2}, $\rm{deg}(\textit{d}_{-})\leq\beta.$  From (\ref{z00}), we have (\ref{djeljivost}), so as in 2.a), $y_{1}=\pm r$ and $z_{1}=\mp s$. Hence, $d_{1}=a$. Using (\ref{prva}) and (\ref{ograda4}), we get $\gamma\geq 2\alpha+\beta$, thus by Lemma \ref{Lm2}, $\rm{deg}(\textit{d}_{-})\geq\alpha.$

\textbf{3. c)} For the case\footnote{A $D(1)$-triple $\{X-1,X+1,16X^3-4X\}$, where $d_{-}=4X$, is an example for this case.} $v_{1}=w_{0}=\pm cr\mp st$, similarly as in 3.b), we get  $x_{0}=\pm r$ and $z_{0}=\mp t$. Hence, $d_{0}=b$ and $d_{1}=d_{-}$. Also, we obtain $\alpha+2\beta\leq \gamma\leq 2\alpha+\beta$ from which it follows that $\alpha=\beta$ and then $\gamma=3\alpha$.

\textbf{3. d)} For the case\footnote{Such a $D(1)$-triple is for example the set $\{\frac{1}{5}, \frac{625X^2+250X}{5}, 12500000000X^{10}+26000000000X^9+23070000000X^8+11392000000X^7+3424950000X^6+644520000X^5+75187000X^4+5200000X^3+194525X^2+3300X+15\}$, where $d_{-}=125000000X^8+210000000X^7+144200000X^6+52040000X^5++10562000X^4+1195600X^3+70160X^2+1800X+56/5$.} $v_{1}=w_{1}=\pm cr\mp st$ we use results from 3.b) and 3.c). They lead to $z_0=\mp t$ and $z_1=\mp s$, so $d_0=b$ and $d_1=a$. Conclusion about degrees follows from (\ref{ograda3}) and (\ref{ograda4}). By Lemma \ref{Lm2}, we have $\rm{deg}(\textit{d}_{-})\geq\beta$.\hfill$\square$

\begin{remark}\label{d-} The case 2.a) of Lemma \ref{Lm4} can be described more precisely. By squaring the equation (\ref{a}), we get $c^2-2bc(a^2+u_{-}^2)+b^2=4a(b+c)+4$, i.e.  \begin{equation} \label{rastav}(c-b(a+u_{-})^2)(c-b(a-u_{-})^2)=4a(b+c)+4. \end{equation}By (\ref{d-plus-minus}), $v_{-}=-r$ and $u_-^2=a^2+1$. From (\ref{ccc}), we get \begin{equation}\label{c-pomocna}c=b+2r^{2}(a-u_{-}).\end{equation}From (\ref{ccc}), using (\ref{prva}) and (\ref{d-plus-minus}), we also get  \begin{equation} \label{jednakost_c}c=b(a - u_{-})^2+2(a - u_{-}),\end{equation}i.e. in (\ref{rastav}) the second factor on left hand side is constant. By (\ref{uvw}), $-r=bs-rt$, so \begin{equation}\label{r1}r(t-1)=bs.\end{equation}By (\ref{prva}), $ \rm{gcd}( \textit{b,r})=1$ so $r=ps$, where $p\in\mathbb{R}^{+}$. Since $b<c$ and $\beta=\gamma$, we have $B\leq C$. By comparing the leading coefficients in (\ref{r1}), we get $p=\frac{\sqrt{B}}{\sqrt{C}}$, so $0<p< 1$ (for $p=1$ we would have $b=c$). Using (\ref{prva}), (\ref{c-pomocna}) and (\ref{d-plus-minus}), since $u_-<0$,  we further conclude that $a-u_{-}=\frac{1}{p}$, i.e. $p=-u_{-}-a$. By (\ref{c-pomocna}), \begin{equation}\label{c-pomocnaa}c=b+\frac{2}{p}r^2=b+2ps^2.\end{equation} Also, from (\ref{r1}), we have \begin{equation}\label{r2}t=\frac{b}{p}+1.\end{equation}From (\ref{c-pomocnaa}), using (\ref{r2}), we obtain \begin{equation}\label{t}c\Big(\frac{1}{p}-2a\Big)=t+1.\end{equation}From (\ref{r2}), (\ref{t}) and (\ref{prva}), we obtain that \begin{equation}\label{1}t+1=cp.\end{equation}By (\ref{t}), (\ref{r2}) and  (\ref{1}), the triple from the case 2.a) of Lemma \ref{Lm4} has the form \begin{equation}\label{abc}\{a,b,c\}=\Big\{\frac{1-p^2}{2p},tp-p, \frac{t}{p}+\frac{1}{p}\Big\}.\end{equation}Also, by (\ref{r2}), \begin{equation}\label{c1}\{a,b,c\}=\Big\{\frac{1-p^2}{2p},b, \frac{b}{p^2}+\frac{2}{p}\Big\}.\end{equation} \end{remark}

In the following lemma, we adjust  \cite[Lemma 10]{dif2} to the situation in $\mathbb{R}[X]$.
\begin{lm}\label{Lm-ubaceno}Let $\{a,b,c\}$, where $a<b<c$, be a polynomial $D(1)$-triple with $\beta< \gamma= \alpha+2 \beta.$ Then $\{a,b,d_-,c\}$ has elements
 \begin{equation}\label{d-_=0}\{a,b,a+b\pm 2r, 4r(r\pm a)(b\pm r)\}\ \ or\end{equation}
\begin{equation}\label{d_-=a}\bigg\{\pm\frac{D_{-}-B}{2\sqrt{BD_-}},b,b\frac{D_-}{B}\pm 2\frac{\sqrt{D_-}}{\sqrt{B}}, \pm2b^2\frac{\sqrt{D_-}}{\sqrt{B}}\Big(\frac{D_-}{B}- 1\Big)+2b\Big(3\frac{D_-}{B}- 1\Big)\pm\frac{9D_{-}-B}{2\sqrt{BD_-}}\bigg\},\end{equation}where $D_-$ is the leading coefficient of $d_-$ and the upper combination of the signs in (\ref{d_-=a}) is for the case $b<d_-$ while the lower is for the case $b>d_-$.
\end{lm}

\emph{Proof.} For the triple $\{a,b,c\}$, by Lemma \ref{Lm4} and Lemma \ref{Lm2},  $\rm{deg}(\textit{d}_{-})=\beta$. Hence,  the triple $\{a,b,d_{-}\}$ has the form 1.) or 2.a) from Lemma \ref{Lm4}. Also, by (\ref{ccc}), for that triple $c=e_+$.

If the triple $\{a,b,d_{-}\}$ is regular, by Definition \ref{def1}, $d_{-}=a+b\pm2r$. Similarly as in \cite[Lemma 10]{dif2}, $c=4r(r\pm a)(b\pm r)$ and $s=2r^2\pm 2ar-1$. Also, by (\ref{w}) we conclude $c=2r(w_{-}\mp 1).$

Let the triple $\{a,b,d_{-}\}$ has the form 2.a) from Lemma \ref{Lm4}, described in Remark \ref{d-}. Then $\alpha=0$. We have $a<b$ and $b<d_-$ or $d_-<b$. Let $b<d_-$ and denote $p_1:=\frac{\sqrt{B}}{\sqrt{D_-}}.$ By (\ref{c-pomocnaa}), $d_-= b+\frac{2r^2}{p_1}$. By (\ref{1}), $v_-=1-p_1d_-$ and we also have $r=-p_1u_-$ (here we use the fact that $u_-,v_-<0$). Using that and (\ref{prva}), from (\ref{ccc}), we obtain $c=a-\frac{4r^2}{p_1}v_-$.  From that, by applying (\ref{r2}) for the triple $\{a,b,d_{-}\}$, we obtain $c=a+\frac{4r^2}{p_1}\big(\frac{b}{p_1}+1\big)$. Using (\ref{prva}) and the expression for $a$ from (\ref{abc}), we get (\ref{d_-=a}). Similarly, for $d_-<b$, we denote $p_2:=\frac{\sqrt{D-}}{\sqrt{B}}$ and we obtain $d_-= b-2p_2r^2$, $v_-=-1-\frac{d_-}{p_2}$ and $u_-=-p_2r$. From that, by applying (\ref{1}) to the triple $\{a,b,d_{-}\}$, we get $c=a+4r^{2}p_2(bp_{2}-1)$. Moreover, applying (\ref{uvw}) for that triple, we obtain $s=-av_{-}-ru_{-}$. It implies that \begin{equation}\label{ss}s=\pm\Big(\frac{D_-}{B}-1\Big )b+\frac{3D_{-}-B}{2\sqrt{BD_{-}}}.\end{equation} By (\ref{w}), $c=a+2r\big(w_{-}\mp r\frac{\sqrt{D_-}}{\sqrt{B}}\big).$ We may also notice that $r|u_-$, so by (\ref{uvw}), $r|t$.\hfill$\square$

\

In Lemma \ref{Lm4} we described different types of polynomial $D(1)$-triples. Note that in cases 1.) and 2.a), $\beta=\gamma$ and in all other cases $\beta<\gamma$. In the rest of the paper we distinguish the cases depending on the parity of indices $m$ and $n$ in the recurring sequences $(v_m)$ and $(w_n)$. From  (\ref{osma}) and (\ref{deveta}), by induction, congruence relations from the following lemma hold for $m,n\geq 0$ (see \cite{dif2}). Here we consider congruences in $\mathbb{R}[X]$.

\begin{lm}\label{Lm5}Let the sequences $(v_{m})$ and $(w_{n})$ be given by  (\ref{osma}) and (\ref{deveta}). Then\begin{eqnarray*}
v_{2m}&\equiv& z_{0}\ ({\rm mod}\ {c}),  \ \ v_{2m+1}\equiv sz_{0}+cx_{0}\ ({\rm mod}\ {c}), \\
w_{2n}&\equiv& z_{1}\ ({\rm mod}\ {c}),  \ \ w_{2n+1}\equiv tz_{1}+cy_{1}\ ({\rm mod}\ {c}).
\end{eqnarray*}\end{lm}

In the following lemma, which is \cite[Lemma 3]{dij}, Dujella and the second author described all possible relations between the initial terms $z_{0}$ and $z_{1}$ of the recurring sequences $(v_m)$ and $(w_n)$ in $\mathbb{C}[X]$. Since every $D(1)$-triple in $\mathbb{R}[X]$ is also a $D(1)$-triple in $\mathbb{C}[X]$, we use that result. Later, we examine all those relations and give some additional information about them which hold in $\mathbb{R}[X]$.

\begin{lm}\label{Lm6} \ \ \\
\emph{\textbf{1)}} If $v_{2m}=w_{2n}$, then $z_{0}=z_{1}.$ \\
\emph{\textbf{2)}} If $v_{2m+1}=w_{2n}$, then either $(z_{0},z_{1})=(\pm
1,\pm s)$ or $(z_{0},z_{1})=(\pm s,\pm 1)$ or
$z_{1}=sz_{0}+cx_{0}$ or $z_{1}=sz_{0}-cx_{0}$.\\
\emph{\textbf{3)}} If $v_{2m}=w_{2n+1}$, then either $(z_{0},z_{1})=(\pm
t,\pm 1)$ or
$z_{0}=tz_{1}+cy_{1}$ or $z_{0}=tz_{1}-cy_{1}$.\\
\emph{\textbf{4)}} If $v_{2m+1}=w_{2n+1}$, then either $(z_{0},z_{1})=(\pm
1, \pm cr\pm st)$ or $(z_{0},z_{1})=(\pm cr \pm st,\pm 1)$ or
$sz_{0}+cx_{0}=tz_{1}\pm cy_{1}$ or $sz_{0}-cx_{0}=tz_{1}\pm cy_{1}$.
\end{lm}

Note that in $\mathbb{Z}[X]$ some relations from Lemma \ref{Lm6} are not possible. By \cite[Lemma 5]{dif2}, if the equation $v_m=w_n$ has a solution, then there exists a solution with $m\in\{0,1\}$. Those solutions induces $d<c$ such that $ad+1$, $bd+1$ and $cd+1$ are perfect squares. In $\mathbb{R}[X]$ in that cases, from minimality assumption, it follows that $d=d_-$ (already described in Lemma \ref{Lm4}) or $d=0\neq d_-$ or $d=a\neq d_-$ where $a$ is a constant (this case is not possible in $\mathbb{Z}[X]$). In $\mathbb{R}[X]$ we also can not exclude all other possibilities from Lemma \ref{Lm6}. Hence, in the next lemma we examine all of them. We determine all possible initial terms and corresponding relations between degrees of polynomials in a polynomial $D(1)$-triple $\{a,b,c\}$.

\begin{lm}\label{Lm7}
\begin{enumerate}
    \item [\textbf{1)}] If $v_{2m}=w_{2n}$, then either
\begin{enumerate}
    \item [\textbf{a)}] $z_{0}=z_{1}= \pm1$ or
 \item [\textbf{b)}] $z_{0}=z_{1}= \pm s$ and $\alpha=0$ or
 \item [\textbf{c)}] $z_{0}=z_{1}= \pm cr \mp st$ and $\alpha>0$, $\alpha+\beta\leq\gamma\leq 2\alpha+\beta$.
\end{enumerate}
\item [\textbf{2)}] If $v_{2m+1}=w_{2n}$, then either
\begin{enumerate}
\item [\textbf{a)}] $z_0=\pm 1$, $z_1=\pm s$ and $\gamma\geq 2\alpha+\beta$ or
\item [\textbf{b)}] $z_0=\pm s$, $z_1=\pm 1$ and $\alpha=0$ or
\item [\textbf{c)}] $z_0=\mp t$, $z_{1}=\pm cr\mp st$ and $\alpha=\beta$, $\gamma=3\alpha$.
\end{enumerate}
\item [\textbf{3)}] If $v_{2m}=w_{2n+1}$, then either
\begin{enumerate}
    \item [\textbf{a)}]  $z_0=\pm t$, $z_1=\pm 1$ and $\gamma\geq \alpha+2\beta$ or
 \item [\textbf{b)}] $z_{0}=\pm s$, $z_{1}=\mp s$ and $\alpha=0$, $\beta=\gamma$ or
\item [\textbf{c)}] $z_{0}=\pm cr\mp st$, $z_{1}=\mp s$ and $\alpha\geq 0$, $2\alpha+\beta\leq\gamma\leq \alpha+2\beta$ or
\item [\textbf{d)}] $z_{0}=\pm s$, $z_{1}=\mp 1$ and $\alpha=0$, $\beta=\gamma$.
\end{enumerate}

\item [\textbf{4)}] If $v_{2m+1}=w_{2n+1}$, then either
\begin{enumerate}
\item [\textbf{a)}] $z_0=\pm 1$, $z_1=\mp cr\pm st$ and $\gamma\leq 2\alpha+\beta$ (with special cases:
\begin{enumerate}
\item [\textbf{1)}] $z_0=\pm 1$, $z_1=\mp 1$ and $\alpha\leq \beta=\gamma$ and
\item [\textbf{2)}] $z_0=\pm 1$, $z_1=\mp s$ and $\alpha=0$, $\beta=\gamma$) or
\end{enumerate}

\item [\textbf{b)}] $z_0=\pm cr\mp st$, $z_1=\mp 1$ and $\gamma\leq \alpha+2\beta$ (with special cases 4.a.1) and\\
$z_0=\pm s$, $z_1=\mp 1$ and $\alpha=0$, $\beta=\gamma$) or
\item [\textbf{c)}] $z_0=\pm t$, $z_1=\pm s$ and $\gamma\geq \alpha+2\beta$.
\end{enumerate}

\end{enumerate}
\end{lm}
\emph{Proof.} \textbf{1)} If $v_{2m}=w_{2n}$ for some integers $m,n\geq 0$, then by Lemma \ref{Lm6}, $z_0=z_1$. Thus, $d_0=d_1=d$ and it holds Remark \ref{rem4}. The cases where $d=d_{-}$ are described in details in 1.)-3.a) of Lemma \ref{Lm4}. If $d=0\neq d_{-}$, then $z_{0}=z_{1}= \pm1$. By Remark \ref{rem3}, if $\beta=\gamma$ then $d_{-}=a$ and $\alpha=0$. Otherwise, $\beta<\gamma$. If $d=a\neq d_{-}$, then $z_{0}=z_{1}= \pm s$ and $\alpha=0$, because we cannot have two different non-constant polynomials in a $D(1)$-quadruple. By Remark \ref{rem3}, if $\beta=\gamma$, then $d_{-}=0$. Otherwise, $\beta<\gamma$. Hence, we get the cases 1.a)-1.c).

\textbf{2.)} By Lemma \ref{Lm6}, if $v_{2m+1}=w_{2n}$, we can have $z_0=\pm 1$, $z_1=\pm s$. By (\ref{ograda4}), $\gamma\geq 2\alpha+\beta$. From (\ref{d_0}), (\ref{d_1}) and Remark \ref{rem2}, we conclude $x_0=1$ and $y_1=r$.  By Lemma \ref{Lm5}, $sz_{0}+cx_{0}\equiv z_{1}\ ({\rm mod}\ {c})$, so if the signs of $z_0$ and $z_1$ are different, then $s=0$ or $c|s$. Those cases are not possible, hence we have 2.a), with equal signs $\pm$.

By Lemma \ref{Lm6}, we further have $z_0=\pm s$, $z_1=\pm 1$. By (\ref{ograda3}), $\gamma\geq 3\alpha$. From (\ref{d_0}), (\ref{d_1}) and Remark \ref{rem2}, we have $x_0^2=a^2+1$ and $y_1=1$. Hence, $\alpha=0$ because otherwise the equation $x_0^2=a^2+1$ is not possible. By Lemma \ref{Lm5}, using (\ref{prva}), if the signs of $z_0$ and $z_1$ are different, then $c|2$ which is not possible. Hence, we have the case 2.b), where the signs $\pm$ are equal.

By Lemma \ref{Lm6} and Remark \ref{rem2}, there is also a possibility $v_{1}=w_0$ where we have to consider both signs in $\pm z_0$ and in $\pm z_1$. Using (\ref{cetvrta}), we get \begin{eqnarray}\label{druga_jednadzba}(sz_{0}+cx_{0})(sz_{0}-cx_{0})=s^{2}z_{0}^{2}-c^2x_{0}^{2}=z_{0}^{2}+ac-c^2. \end{eqnarray} Similarly as in the proof of \cite[Lemma 5]{dif2}, by (\ref{ograda4}) and (\ref{druga_jednadzba}), we conclude that $sz_{0}-cx_{0}=z_{1}$, where  $z_0>0$ and $z_1<0$, or $sz_{0}+cx_{0}=z_{1}$, where  $z_0<0$ and $z_1>0$. By Lemma \ref{Lm3}, $d_1=d_{-}$ or $d_1=0\neq d_{-}$ or $d_1=a\neq d_{-}$ and $\alpha=0$. The case where $d_1=d_{-}$ is described in 3.c) of Lemma \ref{Lm4}. From that we obtain 2.c). If $d_1=0\neq d_{-}$, then $sz_0 \pm cx_0=\pm 1$. If $ \alpha= \gamma$, then by Remark \ref{rem3}, $d_{-}=0$, which is a contradiction. Therefore, $ \alpha< \gamma$, and by (\ref{druga_jednadzba}), $\rm{deg}(\textit{sz}_0\mp \textit{cx}_{0})=2 \gamma$. This is not possible because of (\ref{ograda3}). If $d_1=a\neq d_{-}$ and $\alpha=0$, then $sz_0\pm cx_0=\pm s$. Hence, $s(z_0\mp 1)=\mp cx_0$. We conclude that $s|x_0$, which is not possible since $x_0\neq 0$ and because of (\ref{ograda3}).

\textbf{3.)} If $v_{2m}=w_{2n+1}$, then by Lemma \ref{Lm6}, we have the case 3.a) which is completely analogous to 2.a). Also, by Lemma \ref{Lm6} and Remark \ref{rem2}, we have $v_{0}=w_1$, where we have to consider both signs in $\pm z_0$ and in $\pm  z_1$. Using (\ref{peta}), we get \begin{eqnarray}\label{prva_jednadzba}(cy_{1}+tz_{1})(cy_{1}-tz_{1})=c^{2}y_{1}^{2}-t^2z_{1}^{2}=c^2-bc-z_{1}^2. \end{eqnarray} Similarly as in the proof of \cite[Lemma 5]{dif2},  by (\ref{ograda3}) and (\ref{prva_jednadzba}), we conclude that $z_0=tz_{1}-cy_{1}$, where $z_0<0$ and $z_1>0$, or $z_0=tz_{1}+cy_{1}$, where  $z_0>0$ and $z_1<0$. By Lemma \ref{Lm3}, $d_0=d_{-}$ or $d_0=0\neq d_{-}$ or $d_0=a\neq d_{-}$ and $a$ is a constant. The cases where $d_0=d_{-}$ are described in 2.a) and 3.b) of Lemma \ref{Lm4}. This cases are given in 3.b) and 3.c). If $d_0=0\neq d_{-}$, then $tz_1\pm cy_1=\pm 1$. If $ \beta= \gamma$, then by Remark \ref{rem3}, $d_{-}=a$, $\alpha=0$ and $c=b+ 2rs$. Hence, $\rm{deg}(\textit{tz}_1\mp\textit{cy}_{1})=2 \gamma$. If $ \beta< \gamma$, then we also have $\rm{deg}(\textit{tz}_1\mp\textit{cy}_{1})=2 \gamma$. This is not possible because of (\ref{ograda4}). If $d_0=a\neq d_{-}$ and $\alpha=0$, then $tz_1\pm cy_1=\pm s$.  If $\beta=\gamma$, then by Remark  \ref{rem3}, $d_{-}=0$.  Furthermore, by Lemma \ref{Lm4}, $C=B$ and by (\ref{prva_jednadzba}), $\rm{deg}(\textit{tz}_1\mp \textit{cy}_{1})<  \frac{3\gamma}{2},$ which is possible. If $\beta<\gamma$, then by (\ref{prva_jednadzba}), $\rm{deg}(\textit{tz}_1\mp\textit{cy}_{1})=\frac{3\gamma}{2},$ which is possible only if $\beta=\gamma$, because of (\ref{ograda4}). Hence, we obtain a contradiction. Let $\beta=\gamma$. By Lemma \ref{Lm5}, we have $tz_1\equiv \pm s\ ({\rm mod}\ {c})$. Multiplying that by $t$, we furthermore obtain $z_1\equiv \pm st\pm cr\ ({\rm mod}\ {c})$. Since
\begin{equation}\label{stupnjic}
(\pm st-cr)(\pm st+cr)=ac+bc+1-c^{2},
\end{equation}
one of the
polynomials $\pm st\pm cr$ has degree less then $\gamma$ and the
other has degree equal to $\gamma +\frac{\alpha +\beta}{2}$.
Hence, $\mp st\pm cr =z_{1}$ and ${\rm deg}(z_{1})\leq
\gamma -\frac{\alpha +\beta}{2}$. Also, notice that
$cd_{-}+1=z_{1}^{2}$. If ${\rm deg}(z_{1})< \gamma
-\frac{\alpha +\beta}{2}$, then
${\rm deg}(z_{1})<\frac{\gamma}{2}$, so
$z_{1}=\mp 1$. This is the case 3.d). If ${\rm deg}(z_{1})=\gamma -\frac{\alpha
+\beta }{2}$,
then ${\rm deg}(z_{1})= \frac{\gamma}{2}$. Since
$d_{-}=0$, this is not possible.

\textbf{4.)} If $v_{2m+1}=w_{2n+1}$, then by Lemma \ref{Lm6}, we firstly can have $(z_{0},z_{1})=(\pm
1, \pm cr\pm st)$, i.e. $(z_0,z_1)=(\pm 1, \pm w_{\pm})$. By (\ref{ograda4}), we have $(z_0,z_1)=(\pm 1, \pm w_{-})$, so by Lemma \ref{Lm5}, $(z_0,z_1)=(\pm 1, \mp w_{-})$. As it is described in the proof of \cite[Lemma 3]{dij}, we have the case 4.a). Specially, for $d_{-}=0$, by Lemma \ref{Lm4} and Lemma \ref{Lm5}, $z_1=\mp 1$ and $\alpha\leq\beta=\gamma$. For $d_{-}=a$, we have $(z_0,z_1)=(\pm 1, \mp s)$ and $\alpha=0$, $\beta=\gamma$.

By Lemma \ref{Lm6}, we can also have $(z_{0},z_{1})=(\pm cr\pm st, \pm
1)$, i.e. $(z_0,z_1)=(\pm w_{\pm}, \pm 1)$. By (\ref{ograda3}) and Lemma \ref{Lm5}, we have the case 4.b). Specially, for $d_{-}=0$, $z_0=\pm 1$ and $\alpha\leq\beta=\gamma$, which is the case 4.a.1). For $d_{-}=a$, we have\footnote{In $\mathbb{C}[X]$ here appears an irregular polynomial $D(1)$-quadruple $\mathcal{D}_{p}$ (see \cite[Proposition 1]{dij}).} $(z_0,z_1)=(\pm s, \mp 1)$, $\alpha=0$ and $\beta=\gamma$.

By Lemma \ref{Lm6} and Remark \ref{rem2}, we further have $v_{1}=w_1$, where we have to consider both signs in $\pm z_0$ and $\pm z_1$. Using results from cases 2.) and 3.), we conclude that $sz_{0}-cx_{0}=tz_{1}-cy_1$, where  $z_0>0$ and $z_1>0$, or $sz_{0}+cx_{0}=tz_{1}+cy_{1}$, where  $z_0<0$ and $z_1<0$. By Lemma \ref{Lm3}, these equations can lead to $d_{-}$ or to $d_{+}$  if  $\gamma\geq \alpha+2\beta$ or to an irregular $D(1)$-quadruple $\{a,b,c,d\}$, where $d=0\neq d_{-}$ or $d=a\neq d_{-}$ and $a$ is a constant. Cases where we obtain $d_-$ are described in part 3.d) of Lemma \ref{Lm4}. If $sz_0-cx_0=tz_1-cy_1=\pm w_-$, then we have $z_0=t$, $z_1=s$ and by (\ref{ograda3}), $\gamma \geq \alpha+2 \beta$. Similarly, for $sz_0+cx_0=tz_1+cy_1=\pm w_-$, we have $z_0=-t$, $z_1=-s$ and $\gamma \geq \alpha+2 \beta$. Cases where we obtain $d_+$ are analogous. Therefore, we have 4.c). If $d=0\neq d_{-}$ or $d=a\neq d_{-}$ and $a$ is a constant, then $sz_0\pm cx_0=\pm 1$ or $sz_0\pm cx_0=\pm s$, respectively. Both cases are not possible, as we saw in 2.).\hfill$\square$

\section{PROOF OF THE THEOREM \ref{Tm1}}

We are interested to find all extensions of an arbitrary  polynomial $D(1)$-triple $\{a,b,c\}$ to a polynomial $D(1)$-quadruple. This triple can be extend to a regular quadruple by $d_{-}$ and $d_{+}$, whereas for $d_-=0$ and for $d_-=a$ (a constant), we have improper extensions. We can also have improper irregular $D(1)$-quadruples $\{0,a,b,c\}$ and $\{a,a,b,c\}$, where $a$ is a constant. Moreover, we assumed that we have an irregular polynomial $D(1)$-quadruple $\{a,b,c,d'\}$, such that $0<a<b<c<d'$ and $\delta$ is minimal possible.  By Lemma \ref{Lm1}, we reduced the problem of finding these extensions to the problem of existence of a suitable solution of equation $v_{m}=w_{n}$, where $(v_{m})$ and $(w_{n})$ are binary recurrence sequences defined by (\ref{osma}) and (\ref{deveta}), for some initial values $(z_{0},x_{0})$ and $(z_{1},y_{1})$. In Lemma \ref{Lm7} we described all possible initial terms. We will prove that neither of them leads to the extension with such $d'$.

In the proof of the Theorem \ref{Tm1} we use important relations from the following lemma, obtained by considering the sequences $(v_{m})$ and $(w_{n})$, for $m,n\geq 0$, modulo $4c^{2}$  (see \cite[Lemma 6]{dif2}). Here we again consider congruences in $\mathbb{R}[X]$.

\begin{lm}\label{Lm8}Let the sequences $(v_{m})$ and $(w_{n})$ be given by  (\ref{osma}) and (\ref{deveta}). Then, \begin{eqnarray*}
v_{2m}&\equiv&z_{0}+2c(az_{0}m^2+sx_{0}m)\ ({\rm mod}\ {c^2}), \\
v_{2m+1}&\equiv&s z_{0}+c[2asz_{0}m(m+1)+x_{0}(2m+1)]\ ({\rm mod}\ {c^{2}}),\\
w_{2n}&\equiv&z_{1}+2c(bz_{1}n^{2}+ty_{1}n)\ ({\rm mod}\ {c^{2}}), \\
w_{2n+1}&\equiv&t z_{1}+c[2btz_{1}n(n+1)+y_{1}(2n+1)]\ ({\rm mod}\ {c^{2}}).
\end{eqnarray*}\end{lm} We will also use the following result, which follows directly from (\ref{d-plus}).
\begin{lm}\label{Lm-dodano}Let $\{a,b,c\}$ be (an improper or a proper) polynomial $D(1)$-triple for which (\ref{prva}) holds. Then \begin{equation}\label{lm-dodano}2rst\equiv a+b-d_-\ ({\rm mod}\ {c}).\end{equation}
\end{lm}
\bigskip

\textbf{Proof of Theorem \ref{Tm1}.} Case \textbf{1.a)} $v_{2m}=w_{2n}$, $z_{0}=z_{1}=\pm 1$.\\

By (\ref{d_0}), (\ref{d_1}), Lemma \ref{Lm1} and Remark \ref{rem2}, we have $d_0=d_1=0$, $x_0=1$ and $y_1$=1. By Lemma \ref{Lm8},
\begin{equation}\label{kon1}\pm am^2+sm\equiv \pm bn^2+tn\ ({\rm mod}\ {c}).\end{equation}For $0\in\{m,n\}$ we obtain an improper $D(1)$-quadruple $\{0,a,b,c\}$, which can be regular or irregular. Hence, we assume that $m,n\neq 0$. Similarly as in \cite{dif2}, by (\ref{deseta}) and (\ref{dvanaesta}),
${\rm deg}(v_{2m})=\gamma+(2m-1)\frac{\alpha
+\gamma}{2}$,
${\rm deg}(w_{2n})=\gamma+(2n-1)\frac{\beta
+\gamma}{2}$, except for $\alpha\leq\beta=\gamma$, $z_0=-1$ and $c=a+b+2r$, where ${\rm deg}(v_{2m})=\gamma+(2m-1)\frac{\alpha
+\gamma}{2}$, ${\rm deg}(w_{2n})=\frac{\alpha
+\beta}{2}+(2n-1)\frac{\beta
+\gamma}{2}.$ We also have to consider the case where $c=b+2rs$ and $\alpha=0$, which does not exist in \cite{dif2}. We distinguish subcases $\beta<\gamma$ and $\beta=\gamma$. For $\beta<\gamma$, we obtain a contradiction analogously as in \cite{dif2}.

For $\alpha<\beta=\gamma$, by Lemma \ref{Lm4}, we have $d_-=0$ or $d_-=a$ and $\alpha=0$. For $d_-=0$ we obtain $d=d_+=4r(a+r)(b+r)$ completely analogously as in \cite{dif2}. For $d_-=a$, using (\ref{kon1}), (\ref{r2}), (\ref{c-pomocnaa}) and (\ref{prva}), we obtain $\pm am^2+sm\equiv \mp 2pn^2-n\ ({\rm mod}\ {c}).$ Hence,  $\pm am^2+sm= \mp 2pn^2-n$, which is not possible, because on the right hand side we have a constant and on the left hand side a nonconstant polynomial.

For $\alpha=\beta=\gamma$, completely analogously as in \cite{dif2}, we obtain an improper $D(1)$-quadruple $\{0,a,b,c\}$, which can be regular or irregular or we obtain $d=d_+$.\\

Case \textbf{1.b)} $v_{2m}=w_{2n}$, $z_{0}=z_{1}=\pm s$ and $\alpha=0$.\\

By (\ref{d_0}), (\ref{d_1}), Lemma \ref{Lm1} and Remark \ref{rem2}, $d_0=d_1=a$, $x_0^2=a^2+1$ and $y_1=r$. From Lemma \ref{Lm8} we have
\begin{equation}\label{kon2}\pm asm^2+sx_0m\equiv \pm bsn^2+trn\ ({\rm mod}\ {c}).\end{equation}For $0\in\{m,n\}$ we obtain an improper $D(1)$-quadruple $\{a,a,b,c\}$, which can be regular or irregular. Hence, we assume that $m,n\neq 0$.

By multiplying the congruence (\ref{kon2}) by $s$ and by using  (\ref{prva}) and (\ref{lm-dodano}),  \begin{equation}\label{kon3}\pm am^2+x_0m\equiv \pm bn^2+\frac{an}{2}+\frac{bn}{2}-\frac{d_{-}n}{2}\ ({\rm mod}\ {c}).\end{equation}We separate subcases $\beta<\gamma$ and $\beta=\gamma$. Let $\beta<\gamma$. By Lemma \ref{Lm2}, (\ref{kon3}) implies \begin{equation}\label{jed3}d_-=a+b\pm 2bn\mp 2\frac{am^2}{n}-2\frac{x_0m}{n}.\end{equation}If $\rm{deg}(\textit{d}_{-})<\beta$, then $1\pm 2n=0$, which is not possible. Hence, $\rm{deg}(\textit{d}_{-})=\beta$ and by Lemma \ref{Lm2}, we get $\gamma=2\beta$. Lemma \ref{Lm-ubaceno} implies $d_{-}=a+b\pm 2r$ or $d_-=b\frac{D_-}{B}\pm 2\frac{\sqrt{D_-}}{\sqrt{B}}$. For $d_{-}=a+b\pm 2r$, by (\ref{jed3}), $\pm bn\mp r$ is a constant, which is not possible. For $d_-=b\frac{D_-}{B}\pm 2\frac{\sqrt{D_-}}{\sqrt{B}}$, by (\ref{jed3}), we conclude that $\frac{D_-}{B}=1\pm 2n$. Hence, \begin{equation}\label{razlomak}\frac{D_-}{B}=1+ 2n,\end{equation} and  $z_0=z_1=s$. By (\ref{osma}) and (\ref{deveta}), $v_1=s^2+cx_0$ and $w_1=st+cr$, so $\rm{deg}(\textit{v}_{1})=2\beta$ and $\rm{deg}(\textit{w}_{1})=\frac{5\beta}{2}$. By (\ref{deseta}) and (\ref{dvanaesta}), $\rm{deg}(\textit{v}_{2\textit{m}})=2\beta+(2\textit{m}-1)\beta$ and $\rm{deg}(\textit{w}_{2\textit{n}})=\frac{5\beta}{2}+(2\textit{n}-1)\frac{3\beta}{2}$. From $\rm{deg}(\textit{v}_{2\textit{m}})=\rm{deg}(\textit{w}_{2\textit{n}})$, we get $2m=3n$. By inserting that into (\ref{jed3}), we further conclude that $a-\frac{9an}{2}-3x_0=\pm 2\frac{\sqrt{D_-}}{\sqrt{B}}$. Since the left hand side of this equation is $<0$, we conclude that $d_-=b\frac{D_-}{B}- 2\frac{\sqrt{D_-}}{\sqrt{B}}$. In that case $B>D_-$, thus by (\ref{razlomak}), $2n<0$, which is not possible.

Let now $\beta=\gamma$. By Lemma \ref{Lm4}, $c=a+b+2r$ and $d_-=0$ or $c=b+2rs$ and $d_-=a$. Let $c=a+b+2r$. From (\ref{kon3}), we obtain $\pm am^2+x_0m \mp bn^2-\frac{an}{2}-\frac{bn}{2}=k(a+b+2r)$, where $k\in \mathbb{R}$. Hence, $\mp bn^2-\frac{bn}{2}-kb-2kr$ is a constant. From that, by observing degrees, we get $k=0$ and $\mp n^2-\frac{n}{2}=0$, which is not possible.

Let $c=b+2rs$. By (\ref{jednakost_c}), $c=\frac{b}{p^2}+\frac{2}{p}$. By (\ref{kon3}), we get $\pm am^2+x_0m\equiv \pm bn^2+\frac{bn}{2}\ ({\rm mod}\ {c}),$ so $\pm am^2+x_0m\mp bn^2-\frac{bn}{2}=k(\frac{b}{p^2}+\frac{2}{p})$, where $k\in \mathbb{R}$. From that, by comparing degrees of polynomials, we get
\begin{eqnarray}\label{sustav}\mp n^2-\frac{n}{2}-\frac{k}{p^2}&=&0,\\
\nonumber \pm am^2+x_0m-\frac{2k}{p}&=&0.\end{eqnarray}If we have the upper combination of signs in (\ref{sustav}), then from the first equation we get $k<0$ and from the second equation we get $k>0$, which is a contradiction. The lower combination of signs is possible. Then $z_0=z_1=-s$ and $k>0$. By (\ref{osma}) and (\ref{deveta}), we have $v_1=(x_0-a)c-1$ and $w_1=cr-st$, so $\rm{deg}(\textit{v}_{1})=\gamma$ and by (\ref{stupnjic}),  $\rm{deg}(\textit{w}_{1})=\frac{\gamma}{2}$. From (\ref{deseta}) and (\ref{dvanaesta}), we get $\rm{deg}(\textit{v}_{2\textit{m}})=\gamma+(2\textit{m}-1)\frac{\gamma}{2}$ and $\rm{deg}(\textit{w}_{2\textit{n}})=\frac{\gamma}{2}+(2\textit{n}-1)\gamma$. Moreover, $\rm{deg}(\textit{v}_{2\textit{m}})=\rm{deg}(\textit{w}_{2\textit{n}})$ implies $m=2n-1$. Multiplying the first equation in (\ref{sustav}) with $-2p$ and then by adding those equations, we obtain $m=\frac{x_0-np}{a}$. Since $p=x_0-a$, we have $n-1=\frac{x_0}{a}(1-n)$. For $n>1$, we get $x_0=-a$, which is not possible. For $n=1$, we have $m=1$, thus $z=v_2=w_2$. Similarly as in \cite[Proposition 1]{dij} we obtain $d=d_+$. \\

Case \textbf{1.c)} $v_{2m}=w_{2n}$, $z_{0}=z_{1}=\pm cr \mp st$ and $\alpha>0$, $\alpha+\beta\leq\gamma\leq 2\alpha+\beta$.\\

From (\ref{d_0}), (\ref{d_1}), Lemma \ref{Lm1} and Remark \ref{rem2}, we have $d_0=d_1=d_{-}$ and $x_0=rs-at$, $y_1=rt-bs$. For $\beta=\gamma$, completely analogously as in \cite{dif2}, we obtain $d=d_-$ and $d=d_+$. For $\beta<\gamma$, similarly as in  \cite{dif2}, we obtain $d=d_-$ and $d=d_+$, again. \\

Case \textbf{2.a)} $v_{2m+1}=w_{2n}$, $z_{0}=\pm 1$, $ z_{1}=\pm s$ and $\gamma\geq 2\alpha+\beta$.\\

By (\ref{d_0}), (\ref{d_1}), Lemma \ref{Lm1} and Remark \ref{rem2}, we get $d_0=0$, $d_1=a$, $x_0=1$, $y_1=r$ and we have equal signs for $z_0$ and $z_1$. By Lemma \ref{Lm8}, (\ref{lm-dodano}) and (\ref{prva}), we conclude
\begin{equation}\label{kon5}\pm 2am(m+1)+s(2m+1)\equiv \pm 2bn^2+an+bn-d_-n\ ({\rm mod}\ {c}).\end{equation}We distinguish the cases $\beta<\gamma$ and $\beta=\gamma$. Let $\beta<\gamma$. In this case the congruence (\ref{kon5}) becomes an equation. By Lemma \ref{Lm4}, $\textit{d}_{-}\neq 0$ and, by Lemma \ref{Lm2}, $\rm{deg}(\textit{d}_{-})\geq \alpha$. For $n=0$, from (\ref{kon5}), we get $m=-\frac{1}{2}$ which is not possible, so $n>0$. By (\ref{kon5}),  \begin{equation}\label{jed5}d_-=\pm 2bn+a+b\mp 2a\frac{m(m+1)}{n}-s\frac{2m+1}{n}.\end{equation}For $\beta<\frac{\alpha+\gamma}{2}$, from (\ref{jed5}), we get $D_-=-\sqrt{AC}\frac{2m+1}{n}<0$, which is not possible. Hence, $\beta\geq\frac{\alpha+\gamma}{2}$. If $\beta>\frac{\alpha+\gamma}{2}$, then by (\ref{jed5}), $\rm{deg}(\textit{d}_{-})=\beta$. By Lemma \ref{Lm2}, we get $\gamma=\alpha+2\beta,$ which is not possible. Therefore, $\beta=\frac{\alpha+\gamma}{2}$, i.e. \begin{equation}\label{ograda1a}\gamma= 2\beta-\alpha\end{equation} and $\rm{deg}(\textit{d}_{-})=\beta-2\alpha\leq\beta$. Since  $\gamma\geq 2\alpha+\beta$, we also have
\begin{equation}\label{ograda1b}3\alpha\leq\beta.\end{equation}

For $\alpha=0$, $\gamma=2\beta$ and $\rm{deg}(\textit{d}_{-})=\beta$. By (\ref{jed5}), we have \begin{equation}\label{vodeci_koeficijenti}D_-=B(\pm 2n+1)-\sqrt{AC}\frac{2m+1}{n}.\end{equation} Hence, $z_0=1$ and $z_1=s$. For $\alpha>0$, $\rm{deg}(\textit{d}_{-})<\beta$, thus $B(\pm 2n+1)=\sqrt{AC}\frac{2m+1}{n}.$ Again, we must have $z_0=1$ and $z_1=s$. By (\ref{osma}) and (\ref{deveta}), $v_1=s+c$ and $w_1=st+cr$, so $\rm{deg}(\textit{v}_{1})=\gamma$ and $\rm{deg}(\textit{w}_{1})=\gamma+\frac{\alpha+\beta}{2}$. By (\ref{deseta}) and (\ref{dvanaesta}), \begin{equation}\label{s1}\rm{deg}(\textit{v}_{2\textit{m}+1})=\gamma+2\textit{m}\beta\end{equation} and \begin{equation}\label{s2}\rm{deg}(\textit{w}_{2\textit{n}})=\beta+\textit{n}(\beta+\gamma). \end{equation}From (\ref{ograda1a}), (\ref{s1}) and (\ref{s2}), we get \begin{equation}\label{jednakost-stupnjevi}\alpha(-1+n)=\beta(-1+3n-2m).\end{equation}Let $\alpha=0$ and $\gamma=2\beta$. From (\ref{jednakost-stupnjevi}), we get \begin{equation}\label{diofantskaa}2m-3n=-1.\end{equation}Moreover, we have one of the cases from Lemma \ref{Lm-ubaceno}. If $d_-=a+b \pm 2r$, then $s=2r^2 \pm 2ar-1$. Hence, $D_-=B$ and $C=4AB^2$. From (\ref{vodeci_koeficijenti}) using (\ref{diofantskaa}), we obtain $n=3a$. Furthermore, from  (\ref{jed5}), we get $r(\pm 2\pm 6a)=-\frac{2}{3}m(m+1)-3$. Hence, $-\frac{2}{3}m(m+1)=3$, which is not possible. Let $d_-=b\frac{D_-}{B}\pm 2\frac{\sqrt{D_-}}{\sqrt{B}}$. By comparing the leading coefficients in (\ref{jed5}), using (\ref{diofantskaa}) and the equation \begin{equation}\label{cccc}C=4ABD_-,\end{equation}obtained by comparing the leading coefficients in (\ref{ccc}), we obtain \begin{equation}\label{predznaci}(D_--B)(1\pm 3)=2Bn.\end{equation}If in (\ref{predznaci}) we have the sign $-$, then $d_-<b$ and $-D_-=B(n-1)$. Since $D_->0$, we get $n<1$, which is not possible. If in (\ref{predznaci}), we have the sign $+$, then $d_->b$ and $\frac{D_-}{B}=\frac{n}{2}+1$. Using that, (\ref{ss}) and (\ref{diofantskaa}), from (\ref{jed5}), we obtain $(D_- -B)\big(-\frac{m(m+1)}{n}-4\big)=2D_-$. We get $-\frac{m(m+1)}{n}-4>0$, a contradiction.

Let $\alpha>0$. From Lemma \ref{Lm8}, we obtain
\begin{eqnarray*}2s(am(m+1)-bn^2)\equiv 2trn-(2m+1)\ ({\rm mod}\ {c}).\end{eqnarray*}By squaring that, using (\ref{prva}), we get \begin{equation}\label{kon-duje}4(am(m+1)-bn^2)^2\equiv 4r^2n^2-4trn(2m+1)+(2m+1)^2\ ({\rm mod}\ {c}).\end{equation}Using  (\ref{jed5}), we define \begin{equation}\label{g}g:=(2n+1)b-s\frac{2m+1}{n}=d_{-}-a\Big(1-2\frac{m(m+1)}{n}\Big).\end{equation}By observing degrees of polynomials in (\ref{g}), we conclude that $\rm{deg}(\textit{g})\leq \beta-2\alpha.$ If  $\rm{deg}(\textit{g})<\beta-2\alpha$, then $\beta=3\alpha$ and $D_{-}=A\Big(1-2\frac{m(m+1)}{n}\Big).$ Hence, $1-2\frac{m(m+1)}{n}>0$. By (\ref{jednakost-stupnjevi}), $\alpha=\beta\Big(3+2\frac{1-m}{n-1}\Big)<\beta.$ Hence, $n<m$ and we have $1>2(n+1)$, a contradiction. Therefore, $\rm{deg}(\textit{g})=\beta-2\alpha.$ By (\ref{g}), \begin{equation}\label{b}b=s\frac{2m+1}{n(2n+1)}+\frac{g}{2n+1}.\end{equation}Using (\ref{uvw}) and (\ref{b}), we get \begin{equation}\label{rt} rt=s^2\frac{2m+1}{n(2n+1)}+\frac{gs}{2n+1}-v_{-}.\end{equation}Using (\ref{b}) and (\ref{rt}), from (\ref{kon-duje}) we obtain \begin{eqnarray}\label{kon-velika}\nonumber4a^2m^2(m+1)^2&-&8am(m+1)n^2\Big(s\frac{2m+1}{n(2n+1)}+\frac{g}{2n+1}\Big)+\\\nonumber
&+&4n^4\Big(\frac{(2m+1)^2}{n^2(2n+1)^2}+2sg\frac{2m+1}{n(2n+1)^2}+\frac{g^2}{(2n+1)^2}\Big)\equiv 4r^2n^2-\\ &-&4n(2m+1)\Big(\frac{2m+1}{n(2n+1)}+\frac{gs}{2n+1}-v_{-}\Big)+(2m+1)^2\ ({\rm mod}\ {c}).\end{eqnarray}By considering degrees of polynomials on both sides of the congruence (\ref{kon-velika}), using (\ref{ograda1a}) and (\ref{ograda1b}), we conclude that (\ref{kon-velika}) become an equation. If $\beta>3\alpha$, then by considering the leading coefficients in these equation, we get $2n^2=-2n-1$, which is a contradiction. Let $\beta=3\alpha$. Using (\ref{g}), from the equation obtained from (\ref{kon-velika}), we get that $b\big(-8am(m+1)n^2+8\frac{n^4}{2n+1}g-4an^2+4n^{2}g\big)$ is a polynomial of degree $<\beta$. This is possible only if $a(2m(m+1)+1)=g\big(1+\frac{2n^{2}}{2n+1}\big)$.  Since, $2m(m+1)\neq -1$, it follows that $a|g$. Hence, by (\ref{g}),  $d_{-}=\xi a$, where $\xi\in\mathbb{R}$. By (\ref{d-plus-minus}), $\xi a^2+1=u_{-}^2$. Hence, $\big(a-\frac{u_{-}}{\sqrt{\xi}}\big)\big(a+\frac{u_{-}}{\sqrt{\xi}}\big)=-\frac{1}{\xi}$, which is not possible since both factors on the left hand side of this equation can not be constant.

If $\beta=\gamma$, then $\alpha=0$.  By Lemma \ref{Lm4}, $\textit{d}_{-}=0$ or $\textit{d}_{-}=a$. If $\textit{d}_{-}=0$, then, by (\ref{kon5}), $\pm 2am(m+1)+(a+r)(2m+1)\mp 2bn^2-an-bn=k(a+b+2r)$, where $k\in\mathbb{R}$. By comparing degrees on both sides of this equation, we conclude that $\mp 2n^2-n=k$, i.e. $k$ is an integer. Furthermore, $2m+1=2k$, which is not possible. If $\textit{d}_{-}=a$, by (\ref{kon5}), we have $\pm 2am(m+1)+s(2m+1)\mp 2bn^2-bn=k(b+\frac{2ab}{p}+\frac{2}{p})$, where $k\in\mathbb{R}$. By comparing degrees on both sides of this equation, we conclude that $\mp 2n^2-n=k+\frac{2ka}{p}$. Further, $2m+1=0$, which is not possible. \\

Case \textbf{2.b)}  $v_{2m+1}=w_{2n}$, $z_0=\pm s$, $z_1=\pm 1$ and $\alpha=0$.\\

By (\ref{d_0}), (\ref{d_1}), Lemma \ref{Lm1} and Remark \ref{rem2}, we have $d_0=a$, $d_1=0$, $x_0^2=a^2+1$, $y_1=1$ and we have equal signs for $z_0$ and $z_1$.  By Lemma \ref{Lm8} and (\ref{prva}),
\begin{equation}\label{kon6}\pm a\pm 2am(m+1)+\sqrt{a^2+1}(2m+1)\equiv 2(\pm bn^2+tn)\ ({\rm mod}\ {c}).\end{equation} We distinguish subcases $\beta<\gamma$ and $\beta=\gamma$. If $\beta<\gamma$, the congruence (\ref{kon6}) became an equation. The left hand side of that equation is a constant, so $n=0$. Hence, $z=w_0=z_1=\pm 1$, so we obtain an irregular polynomial $D(1)$-quadruple $\{0,a,b,c\}$.

Let $\beta=\gamma$. By Lemma \ref{Lm4}, $\textit{d}_{-}=0$ or $\textit{d}_{-}=a$. If $\textit{d}_{-}=0$, then, by (\ref{kon6}), (\ref{c-reg}) and (\ref{st-reg}), we have $\pm a\pm 2am(m+1)+\sqrt{a^2+1}(2m+1)- 2(\pm bn^2+(b+r)n)=k(a+b+2r)$, where $k\in\mathbb{R}$. By comparing the leading coefficients on both sides of this equation, we obtain $\mp2n^2-2n=k$ and $-n=k$. Therefore, $n=0$ or $\mp 2n-1=0$, which is not possible. For $n=0$, we get a regular polynomial $D(1)$-quadruple $\{0,a,b,c\}$. If $\textit{d}_{-}=a$, then, by (\ref{kon6}) and Remark \ref{d-}, we get $\pm a\pm 2am(m+1)+\sqrt{a^2+1}(2m+1)- 2\big(\pm bn^2+(\frac{b}{p}+1)n\big)=k\big(\frac{b}{p^2}+\frac{2}{p}\big)$, where $k\in\mathbb{R}$. By comparing the leading coefficients on both sides of this equation,  we obtain $\mp 2n^2-\frac{2n}{p}=\frac{k}{p^2}$ and $\pm a\pm 2am^2\pm 2am+2mx_0+x_0-2n=\frac{2k}{p}.$ From that, \begin{equation}\label{p-ovi}(x_0\pm a)(1+2m)\pm2am^2=\mp4n^2p-2n.\end{equation} Since $x_0=-u_-$ and $0<p<1$, for $z_0=s$ and $z_1=1$, from (\ref{p-ovi}) we get $\frac{1}{p}(1+2m)+2am^2=-4n^2p-2n$, where the left hand side is $>1$ and the right hand side is $\leq 0$, which is not possible. Hence, $z_0=-s$ and $z_1=-1$. In that case, from (\ref{p-ovi}), we get $p(1+2m-4n^2)=2am^2-2n$ and then \begin{equation}\label{p-ovi1}1>2m(am-1)+2n(2n-1).\end{equation} By (\ref{osma}) and (\ref{deveta}), we have $v_1=cp-1$ and $w_1=-t+c$, so $\rm{deg}(\textit{v}_{1})=\rm{deg}(\textit{w}_{1})=\gamma$. By (\ref{deseta}), (\ref{dvanaesta}) and $\rm{deg}(\textit{v}_{2\textit{m}+1})=\rm{deg}(\textit{w}_{2\textit{n}})$, and we obtain \begin{equation}\label{min}m+1=2n.\end{equation} From (\ref{min}) and (\ref{p-ovi1}), we get $1>m^2(2a+1)-m$. Since $a>0$, $1>m(m-1)$, which is not possible for $m\geq 2$. For $m=0$, by (\ref{min}), $n=\frac{1}{2}$ which is also not possible. For $m=1$, by (\ref{min}), $n=1$ and we have $v_3=w_2$. By \cite{dij}, this case is not possible. \\

Case \textbf{2.c)}  $v_{2m+1}=w_{2n}$, $z_0=\mp t$, $z_{1}=\pm cr\mp st$ and $\alpha=\beta$, $\gamma=3\alpha$.\\

By (\ref{d_0}), (\ref{d_1}), Lemma \ref{Lm1} and Remark \ref{rem2}, we have $d_0=b$, $d_1=d_-$, $x_0=r$, $y_1=-v_-=rt-bs$. Since $v_1=\mp st+cr$, we have $\rm{deg}(\textit{v}_{1})=\left\{%
\begin{array}{ll}
    \gamma-\frac{\alpha+\beta}{2}, & \hbox{if $z_0=-t$;} \\
    \gamma+\frac{\alpha+\beta}{2}, & \hbox{if $z_0=t$.} \\
\end{array}%
\right.$ Similarly as in \cite[Lemma 8]{dif2}, $\rm{deg}(\textit{w}_{1})=\left\{%
\begin{array}{ll}
    \frac{3\gamma-\alpha}{2}, & \hbox{if $z_1=cr-st$;} \\
    \frac{\alpha+\gamma}{2}, & \hbox{if $z_1=-cr+st$.} \\
\end{array}%
\right.$ If $z_0<0$ and $z_1>0$, by (\ref{deseta}), (\ref{dvanaesta}) and $\rm{deg}(\textit{v}_{2\textit{m}+1})=\rm{deg}(\textit{w}_{2\textit{n}})$, we obtain \begin{equation}\label{slucaj_a}m=n.\end{equation} If $z_0>0$ and $z_1<0$, similarly, we obtain \begin{equation}\label{slucaj_b}m+1=n.\end{equation}

By Lemma \ref{Lm8} and (\ref{slucaj_a}), we conclude $-2astm(m+1)\equiv -2bst(m^2+m)\ ({\rm mod}\ {c}).$ By multiplying this congruence with $st$, we get $-2am(m+1)\equiv -2b(m^2+m)\ ({\rm mod}\ {c}).$ Since $\beta<\gamma$, from that we obtain $2m(m+1)(-a+b)=0$. For $m=0$, we have $z=z_1=\pm w_-$, so $d=d_-$. The cases $m=-1$ and $a=b$ are not possible. Similarly, by Lemma \ref{Lm8} and (\ref{slucaj_b}), we obtain $m(m+1)(a-b)=0$. For $m=0$, we have $z=v_1=cr+st$, hence $d=d_+$. The cases $m=-1$ and $a=b$ are not possible. \\

Case \textbf{3.a)}  $v_{2m}=w_{2n+1}$,  $z_0=\pm t$, $z_1=\pm 1$ and $\gamma\geq \alpha+2\beta$.\\

By (\ref{d_0}), (\ref{d_1}), Lemma \ref{Lm1} and Remark \ref{rem2}, $d_0=b$, $d_1=0$, $x_0=r$, $y_1=1$ and we have equal signs for $z_0$ and $z_1$.  By Lemma \ref{Lm8}, (\ref{lm-dodano}) and (\ref{prva}), we conclude \begin{equation}\label{kon10}\pm 2am^2+ am+bm-d_-m\equiv \pm 2bn(n+1)+t(2n+1)\ ({\rm mod}\ {c}).\end{equation}In this case $\beta\leq\rm{deg}(\textit{d}_{-})<\gamma$. Hence, (\ref{kon10}) became an equation. If $\rm{deg}(\textit{d}_{-})<\frac{\beta+\gamma}{2}$, then, by considering leading coefficients of polynomials on both sides of these equation, we get $2n+1=0$, a contradiction. If $\rm{deg}(\textit{d}_{-})>\frac{\beta+\gamma}{2}$, then $m=0$, so $z=v_0=z_0=\pm t$. Hence, $d=b$, which is not possible. Therefore, $\rm{deg}(\textit{d}_{-})=\frac{\beta+\gamma}{2}$. By considering the leading coefficients in equation obtained from (\ref{kon10}), we get $-D_-m=\sqrt{BC}(2n+1)$. This is not possible, since on the left hand side of this equation we have a real number $\leq 0$ and the right hand side is $>0$. \\

Case \textbf{3.b)}  $v_{2m}=w_{2n+1}$,  $z_{0}=\pm s$, $z_{1}=\mp s$ and $\alpha=0$, $\beta=\gamma$.\\

By (\ref{d_0}), (\ref{d_1}), Lemma \ref{Lm1} and Remark \ref{rem2}, we have $d_0=d_1=a=d_-$, $x_0^2=a^2+1$, $y_1=r$ and we have different signs for $z_0$ and $z_1$. By (\ref{osma}), $v_1=\pm s^2+cx_0=c(\pm a+x_0)\pm 1$, so $\rm{deg}(\textit{v}_{1})=\gamma$. If $z_0>0$ and $z_1<0$, then, by Lemma \ref{Lm7}, $z_0=tz_1+cy_1=-st+cr$. By (\ref{deveta}), we get $w_1=z_0=s$. Hence, $\rm{deg}(\textit{w}_{1})=\frac{\gamma}{2}$. By (\ref{deseta}) and (\ref{dvanaesta}), $\rm{deg}(\textit{v}_{2\textit{m}})=\gamma+(2\textit{m}-1)\frac{\gamma}{2}$ and $\rm{deg}(\textit{w}_{2\textit{n}+1})=\frac{\gamma}{2}+2\textit{n}\gamma$. From $\rm{deg}(\textit{v}_{2\textit{m}})=\rm{deg}(\textit{w}_{2\textit{n}+1})$, we get $m=2n$. If $z_0<0$ and $z_1>0$, then, by Lemma \ref{Lm7}, $z_0=tz_1-cy_1=st-cr$. By (\ref{deveta}), $w_1=st+cr$. Thus, $\rm{deg}(\textit{w}_{1})=\frac{3\gamma}{2}$. By (\ref{deseta}) and (\ref{dvanaesta}), we have $\rm{deg}(\textit{v}_{2\textit{m}})=\gamma+(2\textit{m}-1)\frac{\gamma}{2}$ and $\rm{deg}(\textit{w}_{2\textit{n}+1})=\frac{3\gamma}{2}+2\textit{n}\gamma$. From $\rm{deg}(\textit{v}_{2\textit{m}})=\rm{deg}(\textit{w}_{2\textit{n}+1})$, we get $m=2n+1$.

By Lemma \ref{Lm8}, we conclude
\begin{equation}\label{kon11}\pm s+2c(\pm asm^2+ sx_0 m)\equiv \mp st+c(\mp 2btsn(n+1)+r(2n+1))\ ({\rm mod}\ {4c^2}).\end{equation} We concluded that $\pm s=\pm cr\mp st$. Using that and by dividing (\ref{kon11}) with $c$, we get \begin{equation}\label{kon12}\pm r\pm 2asm^2+ 2sx_0 m\equiv \mp 2btsn(n+1)+r(2n+1)\ ({\rm mod}\ {4c}).\end{equation}From (\ref{r2}) and (\ref{1}), we have \begin{equation}\label{kon27}bt\equiv 2p\ ({\rm mod}\ {c}).\end{equation}Since in this case $r=ps$, where $p\in\mathbb{R}^{+}$, using (\ref{kon27}) from (\ref{kon12}) we get \begin{equation}\label{komentar}\pm p\pm 2am^2+ 2x_0 m= \mp 4pn(n+1)+p(2n+1).\end{equation}If $z_0>0$ and $z_1<0$, from (\ref{komentar}) using $p=x_0-a$, we get $2n(2an+a)=p(-2n^2-3n)$. This is possible only for $(m,n)=(0,0)$. If $z_0<0$ and $z_1>0$, from (\ref{komentar}) using $p=x_0-a$, we get $-2an(2n+1)=p(2n^2+n)$. This is possible only for $(m,n)=(1,0)$. If $z=v_0=z_0=\pm s$, then $d=a=d_-$. If $z=v_2=w_1$, then $z=cr+ st$, so $d=d_+$. \\

Case \textbf{3.c)}  $v_{2m}=w_{2n+1}$, $z_{0}=\pm cr\mp st$, $z_{1}=\mp s$ and $\alpha\geq 0$,  $2\alpha+\beta\leq\gamma\leq \alpha+2\beta$.\\

By (\ref{d_0}), (\ref{d_1}), Lemma \ref{Lm1} and Remark \ref{rem2}, $d_0=d_-$, $x_0=-u=rs-at$, $d_1=a$, $y_1=r$ and we have different signs for $z_0$ and $z_1$. If $z_0>0$ and $z_1<0$, then, by Lemma \ref{Lm7}, $z_0=tz_1+cy_1=-st+cr$. If $z_0<0$ and $z_1>0$, then, by Lemma \ref{Lm7}, $z_0=tz_1-cy_1=st-cr$. Similarly as in \cite[Lemma 8]{dif2}, if $z_0=cr-st$, then,  $\rm{deg}(\textit{v}_{1})=\frac{3\gamma-\beta}{2}$ so $\rm{deg}(\textit{v}_{2\textit{m}})=\frac{3\gamma-\beta}{2}+(2\textit{m}-1)\frac{\alpha+\gamma}{2}$. Also, Similarly as in \cite[Lemma 8]{dif2}, if $z_0=-cr+st$, then $\rm{deg}(\textit{v}_{1})=\frac{\beta+\gamma}{2}$ so $\rm{deg}(\textit{v}_{2\textit{m}})=\frac{\beta+\gamma}{2}+(2\textit{m}-1)\frac{\alpha+\gamma}{2}$. By (\ref{deveta}), if $z_1<0$, then $w_1=cr-st$, $\rm{deg}(\textit{w}_{1})=\gamma-\frac{\alpha+\beta}{2}$ and $\rm{deg}(\textit{w}_{2\textit{n}+1})=\gamma-\frac{\alpha+\beta}{2}+\textit{n}(\beta+\gamma)$. Also, by (\ref{deveta}), if $z_1>0$, then $w_1=cr+st$, $\rm{deg}(\textit{w}_{1})=\gamma+\frac{\alpha+\beta}{2}$ and $\rm{deg}(\textit{w}_{2\textit{n}+1})=\gamma+\frac{\alpha+\beta}{2}+\textit{n}(\beta+\gamma)$. From $\rm{deg}(\textit{v}_{2\textit{m}})=\rm{deg}(\textit{w}_{2\textit{n}+1})$, for $z_0>0$ and $z_1<0$, we get \begin{equation}\label{stupnjeviA}m(\alpha+\gamma)=n(\beta+\gamma).\end{equation} Also, for $z_0<0$ and $z_1>0$, we get \begin{equation}\label{stupnjeviB}(m-1)(\alpha+\gamma)=n(\beta+\gamma).\end{equation}

By Lemma \ref{Lm8}, similarly as in previous cases, we obtain
\begin{equation}\label{kon13}\mp 2(am (m\pm 1)-bn(n+1))\equiv 2rst\Big(n-m +\Big\{\begin{array}{c}
  0 \\
  1 \\
\end{array}\Big\}\Big)\ ({\rm mod}\ {c}),\end{equation} where the upper case corresponds to (\ref{stupnjeviA}) and the lower corresponds to (\ref{stupnjeviB}).

Assume first that $\gamma=\alpha+2\beta$. By Lemma \ref{Lm-ubaceno}, the triple $\{a,b,d_{-}\}$ has the form 1.) or 2.a) from Lemma \ref{Lm4}. Let $d_{-}=a+b\pm2r$. By (\ref{lm-dodano}), from (\ref{kon13}), we get \begin{equation}\label{kon14}\mp 2(am (m\pm 1)-bn(n+1))=\pm 2r\Big(n-m +\Big\{\begin{array}{c}
  0 \\
  1 \\
\end{array}\Big\}\Big).\end{equation} If $\alpha<\beta$, from (\ref{kon14}), we get $n(n+1)=0$. The case $n=-1$ is not possible and for $n=0$, by (\ref{stupnjeviA}) and (\ref{stupnjeviB}), we get $m=0$ and $m=1$, respectively. For $v_0=w_1$, we have $z=z_0=\pm cr\mp st$ and $d=d_-$. For $v_2=w_1$, by \cite{dij}, we get $d=d_\pm$. If $\alpha=\beta$, then $\gamma=3\alpha$. From (\ref{stupnjeviA}) and (\ref{kon14}), we have $m=n$ and then $m(m+1)(-2a+2b)=0$. Since $m\geq 0$ and $a\neq b$, we get $m=n=0$. Hence, again $z=v_0=w_1$.  From (\ref{stupnjeviB}) and (\ref{kon14}), we have $m=n+1$ and then $n(n+1)(2a-2b)=0$. Since $n\geq 0$ and $a\neq b$, we get $n=0$, $m=1$. Hence, again $z=v_2=w_1$.

If the triple $\{a,b,d_{-}\}$ has the form 2.a) from Lemma \ref{Lm4}, then  $\alpha=0$, $\gamma=2\beta$. By Lemma \ref{Lm-ubaceno}, we have $d_-=b\pm 2r^2\frac{\sqrt{D_-}}{\sqrt{B}}$. Using (\ref{lm-dodano}) and (\ref{prva}), from (\ref{kon13}), we get \begin{equation}\label{kon15}\mp 2(am (m\pm 1)-bn(n+1))= \bigg(a\mp 2(ab+1)\frac{\sqrt{D_-}}{\sqrt{B}}\bigg)\Big(n-m +\Big\{\begin{array}{c}
  0 \\
  1 \\
\end{array}\Big\}\Big).\end{equation} By comparing degrees of polynomials on both sides of the equation (\ref{kon15}), we obtain $a\Big(\mp 2m(m\pm 1)-n+m-\Big\{\begin{array}{c}
  0 \\
  1 \\
\end{array}\Big\}\Big)=\mp 2\frac{\sqrt{D_-}}{\sqrt{B}}\Big(n-m +\Big\{\begin{array}{c}
  0 \\
  1 \\
\end{array}\Big\}\Big)$ and ¸$\pm 2n(n+1)=\mp 2a\frac{\sqrt{D_-}}{\sqrt{B}}\Big(n-m +\Big\{\begin{array}{c}
  0 \\
  1 \\
\end{array}\Big\}\Big)$. From that, it follows \begin{equation}\label{kon16}a^2\Big(\mp 2m (m\pm 1)-n+m-\Big\{\begin{array}{c}
  0 \\
  1 \\
\end{array}\Big\}\Big) \mp 2n(n+1)=0.\end{equation} From (\ref{stupnjeviA}) and (\ref{kon16}), we get $2m=3n$ and then $m=n=0$ (which we already had) or $n>0$ and $a^2= \frac{4(n+1)}{-9n-5}$, a contradiction. From (\ref{stupnjeviB}) and (\ref{kon16}), we get $2m=3n+2$ so $(m,n)=(1,0)$ (which we already had) or $n>0$ and $a^2= \frac{-4(n+1)}{9n+7}$, a contradiction.

Let $\gamma<\alpha+2\beta$. Then, $\alpha<\beta$ and by Lemma \ref{Lm2} and Lemma \ref{Lm4}, we have $\rm{deg}(\textit{d}_{-})<\beta$. Using (\ref{lm-dodano}), from (\ref{kon13}), we get \begin{equation}\label{kon17}\mp 2(am (m\pm 1)-bn(n+1))= (a+b-d_-)\Big(n-m +\Big\{\begin{array}{c}
  0 \\
  1 \\
\end{array}\Big\}\Big).\end{equation} By comparing degrees of polynomials on both sides of the equation (\ref{kon17}), we obtain \begin{equation}\label{kon18+}\pm 2n(n+1)= n-m +\Big\{\begin{array}{c}
  0 \\
  1 \\
\end{array}\Big\}.\end{equation} The first case of (\ref{kon18+}), $m=n(-2n-1)$, is possible only for $m=n=0$. This leads to $d=d_-$. In the second case of (\ref{kon18+}), $m-1=n(3+2n)$. By (\ref{stupnjeviB}), we get $(m,n)=(1,0)$ (which leads to $d=d_{\pm}$) or $3+2n=\frac{\beta+\gamma}{\alpha+\gamma}.$ Since $\beta<\gamma$ or $\alpha=0$ and $\beta=\gamma$, we have $1<3+2n<2$ or $3+2n=2$, respectively, a contradiction in both cases.\\

Case \textbf{3.d)}  $v_{2m}=w_{2n+1}$, $z_{0}=\pm s$, $z_{1}=\mp 1$ and $\alpha=0$, $\beta=\gamma$.\\

By (\ref{d_0}), (\ref{d_1}), Lemma \ref{Lm1} and Remark \ref{rem2}, we have $d_0=a$, $d_1=0$, $x_0^2=a^2+1$, $y_1=1$ and we have different signs for $z_0$ and $z_1$. From the proof of Lemma \ref{Lm7}, we know that $d_-=0$, so $c=a+b+2r$. Using that and (\ref{st-reg}), by Lemma \ref{Lm8}, we obtain  \begin{equation}\label{kon18-a}\pm 1+2s(\pm am^2+ x_0m)\equiv \mp 2btn(n+1)+2n+1\ ({\rm mod}\ {c}).\end{equation} Since, by (\ref{st-reg}), in this case $st\equiv -1 ({\rm mod}\ {c})$, by multiplying (\ref{kon18-a}) by $t$, we obtain  \begin{equation}\label{kon18-b}\pm (b+r)\mp 2am^2-2x_0m\pm 2bn(n+1)-(2n+1)(b+r)=k(a+b+2r),\end{equation} where $k\in\mathbb{R}$. By comparing degrees of polynomials on both sides of the equation (\ref{kon18-b}), we get \begin{eqnarray}\label{usporedbica}\pm 1\pm 2n(n+1)-(2n+1)&=&k,\nonumber\\
\pm 1-(2n+1)&=&2k,\\
\mp 2am^2-2x_0m&=&ka.\nonumber\end{eqnarray} From the first two equations in (\ref{usporedbica}), we get $\pm 1-2n-1=\mp 4n(n+1)$, where only the upper combination of signs is possible, i.e. we have $-n=2n^2.$ Hence, $n=0$ and $k=0$. Using that, from the third equation in (\ref{usporedbica}), we get $am^2+x_0m=0$. Since $am+x_0>0$, we have $m=0$. Therefore, $z=v_0=z_0=\pm s$ and we obtain $d=a$, i.e. $\{a,a,b,c\}$ is an irregular polynomial $D(1)$-quadruple.\\

Case \textbf{4.a)}  $v_{2m+1}=w_{2n+1}$, $z_0=\pm 1$, $z_1=\mp cr\pm st$ and $\gamma\leq 2\alpha+\beta$.\\

By (\ref{d_0}), (\ref{d_1}), Lemma \ref{Lm1} and Remark \ref{rem2}, $d_0=0$, $x_0=1$, $d_1=d_-$, $y_1=rt-bs$ and we have different signs for $z_0$ and $z_1$.  If $z_0>0$ and $z_1<0$, then, by (\ref{osma}), $v_1=s+c$, so $\rm{deg}(\textit{v}_{1})=\gamma$. Also, by (\ref{deveta}), similarly as in \cite[Lemma 8]{dif2}, we get $\rm{deg}(\textit{w}_{1})=\frac{\alpha+\gamma}{2}$. If $z_0<0$ and $z_1>0$, then by (\ref{osma}), we have $v_1=-s+c$. Since $a<c$, $\rm{deg}(\textit{v}_{1})=\gamma$. By (\ref{deveta}), similarly as in \cite[Lemma 8]{dif2}, we get $\rm{deg}(\textit{w}_{1})=\frac{3\gamma-\alpha}{2}$. From (\ref{deseta}), (\ref{dvanaesta}) and $\rm{deg}(\textit{v}_{2\textit{m}+1})=\rm{deg}(\textit{w}_{2\textit{n}+1})$, for $z_0>0$ and $z_1<0$, we get
\begin{equation}\label{stupnjeviA1}m(\alpha+\gamma)=\frac{\alpha-\gamma}{2}+n(\beta+\gamma).\end{equation} Also, for $z_0<0$ and $z_1>0$, we get \begin{equation}\label{stupnjeviB1}m(\alpha+\gamma)=\frac{\gamma-\alpha}{2}+n(\beta+\gamma).\end{equation}

By Lemma \ref{Lm8}, (\ref{lm-dodano}) and (\ref{prva}), similarly as in previous cases, we obtain
\begin{equation}\label{kon18}\pm 2am(m+1)+ s(2m+1)\equiv \Big(\frac{a}{2}+\frac{b}{2}-\frac{d_-}{2}\Big)\Big(2n+ \Big\{\begin{array}{c}
  0 \\
  2 \\
\end{array}\Big\}\Big)+b\Big(\pm 2n^2+\Big\{\begin{array}{c}
  0 \\
  -4n-2 \\
\end{array}\Big\}\Big)\ ({\rm mod}\ {c}).\end{equation} We distinguish the cases $\beta<\gamma$ and $\beta=\gamma$. Let $\beta<\gamma$. We have $\alpha>0$, since otherwise $\gamma\leq \beta$. By Lemma \ref{Lm2} and Lemma \ref{Lm4}, we have $0\leq\rm{deg}(\textit{d}_{-})\leq \alpha$. Hence, in (\ref{kon18}) we have an equation. For $\frac{\alpha+\gamma}{2}>\beta$, by comparing degrees in these equation, we get $2m+1=0$, which is not possible. Hence, $\frac{\alpha+\gamma}{2}\leq\beta$, so $\alpha<\beta$. For  $\frac{\alpha+\gamma}{2}<\beta$, by comparing degrees in equation obtained from (\ref{kon18}), for $z_0>0$ and $z_1<0$, we get $n(1+2n)=0$, i.e. $n=0$. By (\ref{stupnjeviA1}), this is not possible. For  $z_0<0$ and $z_1>0$, we get $2n^2+3n+1=0$, i.e. $n=-\frac{1}{2}$ or $n=-1$, a contradiction in both cases. We are left with the possibility that $\frac{\alpha+\gamma}{2}=\beta$, i.e. \begin{equation}\label{gamica}\gamma=2\beta-\alpha\end{equation}and $\rm{deg}(\textit{d}_{-})=\beta-2\alpha$. From that, we get \begin{equation}\label{alpha-beta}2\alpha\leq\beta\leq 3\alpha.\end{equation} For $z_0>0$ and $z_1<0$, from (\ref{kon18}), we obtain the equation (\ref{jed5}) (for $z_0=1$ and $z_1=s$). By comparing degrees in that equation, we obtain \begin{equation}\label{BC}Bn(2n+1)=\sqrt{AC}(2m+1).\end{equation}Since $m$ is a nonnegative integer, we conclude that $n>0$. By (\ref{stupnjeviA1}) and (\ref{gamica}), we get (\ref{jednakost-stupnjevi}). From Lemma \ref{Lm8}, similarly as in the case 2.a), we obtain the congruence (\ref{kon-duje}). Again, we use a polynomial $g$ defined by (\ref{g}), by which from (\ref{kon-duje}) we get the congruence (\ref{kon-velika}). By considering degrees of polynomials which appear in (\ref{g}), similarly as for the case 2.a), we conclude that $\rm{deg}(\textit{g})=\alpha$. By considering degrees of polynomials which appear in (\ref{kon-velika}), we conclude that if $\beta>2\alpha$ then (\ref{kon-velika}) become an equation. By considering leading coefficients in that equation,  we get a contradiction similarly as in the case 2.a). By (\ref{alpha-beta}), we are left with the possibility that $\beta=2\alpha$ and $d_{-}$ is a nonzero constant. In this case, we can apply Lemma \ref{Lm-ubaceno} to a $D(1)$-triple $\{d_{-},a,b\}$ and we have to observe two possibilities.

The first possibility, by Lemma \ref{Lm-ubaceno}, is that we have \begin{equation}\label{b1}b=4u_{-}(u_{-}\pm d_{-})(a\pm u_{-}).\end{equation}By considering leading coefficients of polynomials on both sides of the equation (\ref{b1}), we conclude  that $B=4A^2D_{-}$. Using that and (\ref{cccc}), from (\ref{BC}), we obtain \begin{equation}\label{mn}n(2n+1)=2m+1.\end{equation}For $\beta=2\alpha$, from (\ref{jednakost-stupnjevi}), we get \begin{equation}\label{diofantskaaa}5n-4m=1.\end{equation}All nonnegative solutions $(m,n)$ of the equation (\ref{diofantskaaa}) are given with \begin{equation}\label{min}m=1+5t,\ \ n=1+4t,\end{equation} where $t\in\mathbb{N}_{0}$. By inserting (\ref{min}) into (\ref{mn}), we obtain $t(16t+5)=0$. Hence, $t=0$ i.e. $(m,n)=(1,1)$ or $t=-\frac{5}{16}$, which is not possible. For $(m,n)=(1,1)$ we have $v_3=w_3$.  In this case, by (\ref{osma}) and (\ref{deveta}) \begin{equation}\label{33}sz_0+c(4asz_0+3x_0)+4ac^2x_0=tz_1+c(4btz_1+3y_1)+4bc^2y_1.\end{equation} From (\ref{33}), using (\ref{prva}), we get $2(bs+rt)=4as+4s^2-1$. Therefore, \begin{equation}\label{s}2(rt-bs)\equiv -1\ ({\rm mod}\ {s}).\end{equation} Using (\ref{prva}), we conclude that\begin{equation}\label{bsrt}(bs+rt)(bs-rt)=b^2-ab-bc-1.\end{equation} From (\ref{bsrt}), we furthermore conclude that $\rm{deg}(\textit{bs}-\textit{rt})=\frac{\gamma-\alpha}{2}$, thus, by (\ref{s}), $2(rt-bs)=-1$, i.e. $bs-rt=v_->0$, which is a contradiction.

The second possibility, by Lemma \ref{Lm-ubaceno}, is that we have \begin{equation}\label{d-p}d_{-}=\frac{1}{2}\Big(\frac{1}{p}-p\Big)\end{equation}and \begin{equation}\label{b2-1}b=d_{-}+4u_{-}^2p(ap-1)\end{equation}if $e_{-}<a$ or \begin{equation}\label{b2-2}b=d_{-}+4\frac{u_{-}^2}{p}\Big(\frac{a}{p}+1\Big)\end{equation}if $e_{-}>a$, where $0<p<1$ and $e_{-}$ is obtained by (\ref{d-plus}) for the triple $\{d_{-},a,b\}$. By considering degrees of polynomials on both sides of the equation (\ref{b2-1}), we get $B=4A^2D_{-}p^2$. Using that and (\ref{cccc}), from (\ref{BC}), we obtain \begin{equation}\label{p1}p=\frac{2m+1}{n(2n+1)}.\end{equation}Analogously, from  (\ref{b2-2}), we get $B=\frac{4A^2D_{-}}{p^2}$ and then \begin{equation}\label{p2}p=\frac{n(2n+1)}{2m+1}.\end{equation}From (\ref{ccc}), using (\ref{d-plus-minus}) and (\ref{uvw}) for the triple $\{d_{-},a,b\}$, we obtain \begin{equation}\label{c-treca}c=-a+b+d_{-}-2sv_{-}.\end{equation}By considering degrees of polynomials which appear in (\ref{kon-velika}), using (\ref{c-treca}), we get \begin{eqnarray*}4a^2m^2(m+1)^2˛-8am(m+1)n^2\Big(s\frac{2m+1}{n(2n+1)}+\frac{g}{2n+1}\Big)+\end{eqnarray*} \begin{equation}\label{jed-velika}+4n^4\Big(\frac{(2m+1)^2}{n^2(2n+1)^2}+2sg\frac{2m+1}{n(2n+1)^2}+\frac{g^2}{(2n+1)^2}\Big)- 4r^2n^2+\end{equation} \begin{eqnarray*}+4n(2m+1)\Big(\frac{2m+1}{n(2n+1)}+\frac{gs}{2n+1}-v_{-}\Big)-(2m+1)^2=l(-a+b+d_{-}-2sv_{-}),\end{eqnarray*}where $l\in\mathbb{R}$. A polynomial $v_{-}$ plays for the triple $\{d_{-},a,b\}$ the same role as polynomial $s$ does for the triple $\{a,b,c\}$. Therefore, by (\ref{ss}) (where $s>0$ and $v_{-}<0$), for (\ref{b2-1}) and (\ref{p1}) we have \begin{equation}\label{v1}v_{-}=(p^2-1)a-\frac{1}{2}\Big (3p-\frac{1}{p}\Big).\end{equation}Analogously, for (\ref{b2-2}) and (\ref{p2}) we have \begin{equation}\label{v2}v_{-}=\Big(1-\frac{1}{p^2}\Big)a-\frac{1}{2}\Big (3\frac{1}{p}-p\Big).\end{equation}From (\ref{g}), we get \begin{equation}\label{gg}g=d_{-}-a\Big(1-2\frac{m(m+1)}{n}\Big).\end{equation}Also, from (\ref{g}), for (\ref{b2-1}) and (\ref{p1}) we have \begin{equation}\label{ss1}s=\frac{b}{p}-\frac{gn}{2m+1},\end{equation}and for (\ref{b2-2}) and (\ref{p2}) we have \begin{equation}\label{ss2}s=pb-\frac{gn}{2m+1}.\end{equation}By inserting (\ref{min}), (\ref{p1}), (\ref{d-p}), (\ref{uvw}), (\ref{b2-1}), (\ref{gg}), (\ref{v1}) and (\ref{ss1}) into the equation (\ref{jed-velika}), we get the expression of the form $c_3a^3+c_2a^2+c_1a+c_0$, where $c_i$ are rational expressions with unknowns $l$ and $t$, for $i=0,...,3$.  By solving the system $c_3=0$, $c_2=0$ in unknowns $l$ and $t$, the only integer $t$ we obtain is $t=0$. But, for $t=0$ the coefficients $c_1$ and $c_0$ can not both be equal to $0$. Completely analogously we conclude if we insert  (\ref{min}), (\ref{p2}), (\ref{d-p}), (\ref{uvw}), (\ref{b2-2}), (\ref{gg}), (\ref{v2}) and (\ref{ss2}) into the equation (\ref{jed-velika}).

For $z_0<0$ and $z_1>0$, by comparing degrees in the equation obtained from (\ref{kon18}), we get $\sqrt{AC}(2m+1)=B(-2n^2-3n-1)<0$, which is not possible.

If $\beta=\gamma$, for $d_-=0$ we have $z_1=\mp 1$, $\alpha\leq\beta$, and for $d_{-}=a$ we have $z_1= \mp s$, $\alpha=0$. Let $d_-=0$. Then, $d_1=0$ and $y_1=1$. By Lemma \ref{Lm8}, using (\ref{c-reg}), (\ref{st-reg}) and the fact that $c=s+t$, similarly as in previous cases, we obtain
\begin{equation}\label{kon19}\pm 1 \pm 2asm(m+1)+ 2m\equiv \pm 2bsn(n+1)+2n\ ({\rm mod}\ {c}).\end{equation} Multiplying (\ref{kon19}) by $s$ and by using (\ref{prva}), we get \begin{equation}\label{kon20}(\pm 1 +2m-2n)s \pm 2am(m+1)\mp 2bn(n+1)=k(a+b+2r),\end{equation} where $k\in\mathbb{R}$. If $\alpha<\beta$, by comparing degrees in (\ref{kon20}), we first get $k=\mp 2n(n+1)$. Further, $\pm 1+2m-2n=\mp 4n(n+1)$. If $n\geq 1$, then $2|(\pm 1)$, which is not possible. For $n=0$, we have $\pm 1+2m=0$, which is also not possible. Therefore,  $\alpha=\beta=\gamma$ and $C=A+B+2\sqrt{AB}$. By (\ref{osma}) and (\ref{deveta}), we have $\rm{deg}(\textit{v}_{1})=\rm{deg}(\textit{w}_{1})=\gamma$. By (\ref{deseta}), (\ref{dvanaesta}) and $\rm{deg}(\textit{v}_{2\textit{m}+1})=\rm{deg}(\textit{w}_{2\textit{n}+1})$, we obtain $m=n$. Using that and the fact that $b-a=c-2s$, from (\ref{kon20}), we get \begin{equation}\label{kon21}s(\pm 1 \pm 4n(n+1))\equiv 0\ ({\rm mod}\ {c}).\end{equation} Since $\rm{gcd}(\textit{c,s})=1$, we have $\pm 1 \pm 4n(n+1)=0$, which is a contradiction.

If $d_-=a$, then $d_1=a$ and $y_1=r$. By Lemma \ref{Lm8}, (\ref{prva}), (\ref{c-pomocnaa}), (\ref{r2}) and the fact that $r=ps$, where $p\in\mathbb{R}$ and $0<p<1$, similarly as in previous cases, we obtain
\begin{equation}\label{kon22}\pm 2am(m+1)+ s(2m+1)\equiv \mp \frac{1}{p}\pm 2a \mp 4n(n+1)p+p(2n+1)\ ({\rm mod}\ {c}).\end{equation} Since degrees of polynomials on both sides of the congruence (\ref{kon22}) are less than $\gamma$, we get an equation. By considering degrees of polynomials on both sides of that equation, we get $2m+1=0$, which is not possible. \\

Case \textbf{4.b)}  $v_{2m+1}=w_{2n+1}$, $z_0=\pm cr\mp st$, $z_1=\mp 1$ and $\gamma\leq \alpha+2\beta$.\\

By (\ref{d_0}), (\ref{d_1}), Lemma \ref{Lm1} and Remark \ref{rem2}, $d_0=d_-$, $x_0=rs-at$, $d_1=0$, $y_1=1$ and we have different signs for $z_0$ and $z_1$. By Lemma \ref{Lm8}, (\ref{lm-dodano}) and (\ref{prva}), we have
\begin{equation}\label{kon22-1}\mp 2am (m+ 1)+\Big(-\frac{a}{2}+\frac{b}{2}-\frac{d_-}{2}\Big)\Big(2m+\Big\{\begin{array}{c}
  2 \\
  0 \\
\end{array}\Big\}\Big)\equiv  \mp 2bn(n+1)+t(2n+1)\ ({\rm mod}\ {c}),\end{equation}where the upper case is for $z_0>0$, $z_1<0$ and the lower is for $z_0<0$, $z_1>0$.

If $\beta<\gamma$,  from (\ref{kon22-1}) we obtain the equation. Since $\rm{deg}(\textit{d}_{-})\leq \beta$ and $\rm{deg}(\textit{t})> \beta$, by comparing degrees in these equation, we get $2n+1=0$, which is not possible.

Let $\beta=\gamma$. For $d_-=0$, we have $z_0=\pm 1$, $z_1=\mp 1$ and $\alpha\leq\beta$. By (\ref{d_0}), (\ref{d_1}) and Lemma \ref{Lm1}, we get $d_0=0$, $x_0=1$, $d_1=0$, $y_1=1$. By Lemma \ref{Lm8} and (\ref{prva}), using (\ref{c-reg}), (\ref{st-reg}) and $st\equiv  -1\ ({\rm mod}\ {c})$, similarly as in previous cases, we obtain
\begin{equation}\label{kon23}b(\pm 1+2m\pm 2n(n+1)-2n)+r(\pm 1+2m-2n)+a(\mp 2m(m+1))=k(a+b+2r),\end{equation} where $k\in\mathbb{R}$. If $\alpha<\beta$, by comparing the leading coefficients in (\ref{kon23}), we first obtain $\pm 1+2m\pm 2n(n+1)-2n=k$. Then, $\pm 1+2m-2n=2k$, so $k=\mp 2n(n+1)$. Furthermore, we get $k=\mp 2m(m+1)$. Therefore, $n=m$ or $n=-m-1$. For $m=n$, we get $k=\pm \frac{1}{2}$ and then, $\mp 2m(m+1)= \pm \frac{1}{2}$. This is not possible. The case where $n=-m-1<0$ is also not possible. Hence, $\alpha=\beta=\gamma$. Since $a<b<c$, from (\ref{osma}) and (\ref{deveta}), we conclude that $\rm{deg}(\textit{v}_{1})=\rm{deg}(\textit{w}_{1})=\gamma$. By (\ref{deseta}), (\ref{dvanaesta}) and $\rm{deg}(\textit{v}_{2\textit{m}+1})=\rm{deg}(\textit{w}_{2\textit{n}+1})$, we obtain $m=n$. Using that and the fact that $b-a=c-2s$, from (\ref{kon23}), we get (\ref{kon21}) which is not possible as in the case 4.a).

For $d_-=a$,  we have $z_0=\pm s$, $z_1=\mp 1$, $\alpha=0$ and $\beta=\gamma$. By (\ref{d_0}), (\ref{d_1}), Lemma \ref{Lm1} and Remark \ref{d-}, we have $d_0=a$, $x_0=a+p$, where $p\in\mathbb{R}$ and $0<p<1$, and $d_1=0$, $y_1=1$. By Lemma \ref{Lm8}, \begin{eqnarray*}\pm s^2+c(\pm 2as^2 m(m+1)+(a+p)(2m+1))\equiv \mp t+c(\mp 2btn(n+1)+2n+1)\ ({\rm mod}\ {c^2}).\end{eqnarray*} From (\ref{t}),  $c\Big(\frac{1}{p}-a\Big)=t+s^2.$ Using that, and dividing the obtained congruence by $c$, we get \begin{equation}\label{kon26}\pm \Big(\frac{1}{p}-a\Big)\pm 2as^2 m(m+1)+(a+p)(2m+1)\equiv \mp 2btn(n+1)+2n+1\ ({\rm mod}\ {c}).\end{equation} From (\ref{a}), we obtain $\frac{1}{p}-a=p+a$. Using that, (\ref{kon27}) and (\ref{prva}), by (\ref{kon26}), we conclude \begin{equation}\label{kon28}(p+a)(\pm 1+2m+1)\pm 2am(m+1)=\mp 4pn(n+1)+2n+1.\end{equation} By (\ref{osma}), $v_1=\pm 1+c(\pm a+a+p)$ and by (\ref{deveta}) and (\ref{1}), $w_1=t(\mp 1+\frac{1}{p})+\frac{1}{p}$. By Remark \ref{d-}, $\rm{deg}(\textit{v}_{1})=\rm{deg}(\textit{w}_{1})=\gamma$. From (\ref{deseta}), (\ref{dvanaesta}) and $\rm{deg}(\textit{v}_{2\textit{m}+1})=\rm{deg}(\textit{w}_{2\textit{n}+1})$, we obtain $m=2n$. Using that, for $z_0>0$, from (\ref{kon28}), we get \begin{equation}\label{kon29}(2n+1)(2p+2a+4an-1)=-4pn(n+1).\end{equation} For $(m,n)=(0,0)$, we have a $D(1)$-quadruple $\mathcal{D}_{p}$ from \cite{dij}, whose elements are not from $\mathbb{R}[X]$. For $n>0$, the right hand side in (\ref{kon29}) is $<0$, so we conclude that $4an<1-2x_0<-1$, which is not possible. For $z_0<0$, from (\ref{kon28}), we get \begin{eqnarray*}-4n^2(2a+p)=2n+1,\end{eqnarray*} where the left hand side is $\leq 0$ and the right hand side is $>0$, a contradiction. \\

Case \textbf{4.c)}  $v_{2m+1}=w_{2n+1}$, $z_0=\pm t$, $z_1=\pm s$ and $\gamma\geq \alpha+2\beta$.\\

By (\ref{d_0}), (\ref{d_1}), Lemma \ref{Lm1} and Remark \ref{rem2}, $d_0=b$, $x_0=r$, $d_1=a$, $y_1=r$ and we have equal signs for $z_0$ and $z_1$. Similarly as in \cite{dif2}, we obtain
\begin{equation}\label{kon31}\pm 2astm (m+ 1)+2mr\equiv  \pm 2bstn(n+1)+2nr\ ({\rm mod}\ {c})\end{equation} and \begin{equation}\label{kon32}m(\alpha+\gamma)=n(\beta+\gamma).\end{equation} In this case $\beta<\gamma$. For $\alpha=\beta$, from (\ref{kon31}) and (\ref{kon32}), we obtain $\pm 2m (m+ 1)(a-b)\equiv  0\ ({\rm mod}\ {c}).$ This is possible only for $(m,n)=(0,0)$, where $z=v_1=w_1$, i.e. $d=d_{\pm}$.

Let $\alpha<\beta<\gamma$. By multiplying (\ref{kon31}) with $st$ and by (\ref{lm-dodano}), we get
\begin{equation}\label{kon33}\pm 2am (m+ 1)\mp 2bn(n+1)\equiv  (a+b-d_-)(n-m)\ ({\rm mod}\ {c}).\end{equation} Since $\rm{deg}(\textit{d}_{-})<\gamma$, the congruence (\ref{kon33}) becomes an equation. In this case $\rm{deg}(\textit{d}_{-})\geq \beta$. For $\rm{deg}(\textit{d}_{-})> \beta$, by considering the leading coefficients of polynomials in obtained equation, we get $m=n$, thus, by (\ref{kon32}), $\alpha=\beta$, which is not possible. Therefore, $\gamma=\alpha+2\beta$. We apply Lemma \ref{Lm-ubaceno}. If $d_-=a+b\pm 2r$, then, by considering the leading coefficients of polynomials on both sides of the equation obtained from (\ref{kon33}), we conclude that $\pm n(n+1)=0$. Therefore, $(m,n)=(0,0)$, which leads to $d=d_{\pm}$. If $d_-=b\pm 2r^2\frac{\sqrt{D_-}}{\sqrt{B}}$, similarly, we conclude that $\mp 2n(n+1)=\pm 2a\frac{\sqrt{D_-}}{\sqrt{B}} (n-m)$ and $\pm 2am(m+1)=a\pm 2\frac{\sqrt{D_-}}{\sqrt{B}} (n-m)$, where the signs on the right hand sides of these equations does not depend on the signs on the left hand sides, but are the same in both equations. From that, we obtain $2a^2 m(m+1)=a^2-2n(n+1)$ or $-2a^2 m(m+1)=a^2+2n(n+1)$. In the first case $2n(n+1)=a^2(1-2m(m+1))$, which is not possible. In the second equation, on the left hand side we have a real number $\leq 0$, and on the right hand side we have a positive real number, which is not possible.\hfill$\square$\\

\bigskip

\textbf{Acknowledgement:} The authors were supported by Croatian Science Foundation under the project no. 6422. The second author was also supported by the University of Rijeka research grant no. 13.14.1.2.02. \\

\bigskip
Faculty of Civil Engineering, University of Zagreb,\\ Fra Andrije Ka\v{c}i\'{c}a-Mio\v{s}i\'{c}a 26, 10000 Zagreb, Croatia \\
Email: filipin@grad.hr \\
\\
Department of Mathematics, University of
Rijeka,\\ Radmile Matej\v{c}i\'{c} 2, 51000 Rijeka, Croatia \\
Email: ajurasic@math.uniri.hr

\end{document}